\documentclass[12pt,a4paper]{article}
\usepackage{mathrsfs}
\usepackage{epsfig, graphicx}
\usepackage{latexsym,amsfonts,amsbsy,amssymb}
\usepackage{amsmath,amsthm}
\usepackage{color}
\usepackage[colorlinks, citecolor=blue]{hyperref}
\makeatletter 

\@addtoreset{equation}{section}
\makeatother 
\allowdisplaybreaks 

\newtheorem{remark}{Remark}[section]

\usepackage{subfigure}

\newcommand{\pa}{\partial}

\newcommand{\Om}{\Omega}

\newcommand{\Lam}{\Lambda}

\newcommand{\rtri}{\nabla}

\newcommand{\R}{{\mathbb{R}}}

\newcommand{\be}{\begin{eqnarray}}
\newcommand{\ee}{\end{eqnarray}}
\newcommand{\beq}{\begin{equation}}
\newcommand{\eeq}{\end{equation}}
\newcommand{\ben}{\begin{eqnarray*}}
\newcommand{\een}{\end{eqnarray*}}
\newcommand{\bsp}{\begin{split}}
\newcommand{\esp}{\end{split}}
\newcommand{\nn}{\nonumber}

\numberwithin{equation}{section}
\numberwithin{table}{section}

\begin{document}
\graphicspath{{figures/}} 

\title{A weighted first-order formulation for solving anisotropic diffusion equations with
deep neural networks}

\author{
Hui Xie\footnote{Institute of Applied Physics and Computational Mathematics,
                       Beijing, 100094, China(xiehui@lsec.cc.ac.cn)},
Chuanlei Zhai\footnote{Institute of Applied Physics and Computational Mathematics,
                       Beijing, 100094, China(zhai\_chuanlei@iapcm.ac.cn)},
Li Liu\footnote{Institute of Applied Physics and Computational Mathematics,
                       Beijing, 100094, China(liu\_li@iapcm.ac.cn)},
Heng Yong\footnote{Institute of Applied Physics and Computational Mathematics,
                       Beijing, 100094, China(yong\_heng@iapcm.ac.cn)}
}
\date{}
\maketitle
\begin{abstract}
In this paper, a new weighted first-order formulation is proposed for solving the
   anisotropic diffusion equations with deep neural networks.
   For many numerical schemes, the accurate approximation of anisotropic
   heat flux is crucial for the overall accuracy.
   In this work, the heat flux is firstly decomposed into two components along
   the two eigenvectors of the diffusion tensor, thus the anisotropic
   heat flux approximation is converted into the approximation of
   two isotropic components.
   Moreover, to handle the possible jump of the diffusion tensor across the interface,
   the weighted first-order formulation is obtained by multiplying this
   first-order formulation by a weighted function. By the decaying property of the
   weighted function, the weighted first-order formulation is always well-defined in
   the pointwise way. Finally, the weighted first-order formulation
   is solved with deep neural network approximation.
   Compared to the neural network approximation with the original
   second-order elliptic formulation,
   the proposed method can significantly improve the accuracy, especially for the
   discontinuous anisotropic diffusion problems.


\end{abstract}

\section{Introduction}
\label{intro}
Anisotropic diffusion equations arise from numerical simulations of many important and
practical applications such as inertial confinement fusion (ICF) \cite{Lan,Lindl}, reservoir
simulation \cite{schneider}, astrophysical fluid dynamics \cite{meyer2012}.
Generally, the modelling of the anisotropic diffusion equation is coupled with
the fluid dynamics. Particularly, in the Lagrangian reference frame, the
computational mesh is moved with the fluid motions. For the traditional discrete schemes,
the accuracy is degenerated due to the the distorted meshes
and the anisotropic diffusion tensors. In addition, the different materials or
local physical field (e.g., magnetic field) in the
problems will result in the high discontinuous diffusion tensor across the
interface. Even with the Eulerian reference frame, the anisotropic and discontinuous
diffusion tensor still exists.

Given the anisotropic diffusion tensor and possible discontinuity across the
interface, how to design the accurate discretization for anisotropic
diffusion problems have been investigated in the last decades.
Finite volume (FV) discretizations of such practical multi-physics problems \cite{yong} are
very popular due to easy implementation for irregular computational grid and local conservation.
As we know, the popular two-point flux approximation (TPFA) has been widely used,
nevertheless, its accuracy will degenerate with the presence of anisotropic diffusion tensor.
To recover the (second-order) accuracy in the anisotropic case, some finite
volume schemes based on more than two points on each side are designed.
The pioneer work was proposed in \cite{kershaw1981}, which is the so-called
Kershaw scheme. After that, the local support operator scheme (or mimetic finite difference
method) \cite{Morel}, nine-point scheme \cite{chen2003,sheng0}, multi-point flux
approximation (MPFA) \cite{Aavatsmark} were introduced for the accuracy
improvement. Also some variants \cite{Edwards,gao2} were proposed
for highly anisotropic diffusion tensors.
In addition, some nonlinear schemes \cite{potier,lipnikov1,terekov,xie2, zhao} were proposed for preserving the non-negativity or the discrete maximum principle. However, these traditional
numerical schemes have some restrictions on the the mesh distortion and
anisotropic diffusion tensor, which limits their application range more or less.

With the
advent of the powerful computer and big data, deep neural networks
have been extensively used in diverse tasks such as image classification,
computer vision, language translation, game intelligence, reservoir parameter inversion \cite{xiong}.
In the context of scientific computing, the idea to use neural
network for solving PDEs has been raised since 1990s \cite{Lee, lagaris1998}.
Contrast to the traditional mesh-based numerical algorithms,
the neural network can be seen as a meshless optimization method.
The training process is based on the random sampling strategy.
Thus it is less sensitive to the dimensionality of the problems.
The high-dimensional PDEs that can not be tackled with traditional
mesh-based numerical algorithms are solved by deep neural network
\cite{Han,Sirignano}.

Recently, Raissi et al. \cite{Raissi} developed the so-called
physics-informed neural network (PINN) for the solution and discovery of PDEs.
the main idea behind the PINN is that the governing equation is used, rather
than the labeled solution, in the
loss function to keep the neural network solution approaching
the strong solution of PDEs. PINNs have been successfully used in
solving a large number of nonlinear PDEs, including Burgers,
Schr\"{o}inger, Navier-Stokes, Allen-Cahn, high-speed flows, etc \cite{Lulu, Wight, mao, Cuomo}.

As for the general second-order elliptic problems, the
convergence of PINN was studied for linear second-order
elliptic and parabolic case firstly in \cite{Shin} and
then in \cite{Jiao} for a quantitative error estimation with respect to
the neural networks' depth, width and sampling numbers.
Solving a class
of second-order boundary-value problems on complex geometries
was studied with the penalty manner \cite{Berg} and penalty-free manner\cite{Shenghl}.
A class of linear and nonlinear elliptic problems with
the unknown diffusion coefficient was considered and learned
by the PINN in \cite{Tartakovsky}. In \cite{He}, the forward and
backward advection-diffusion equations are addressed by PINN, especially for high
P\'{e}clet number problems. A extreme learning machine with the single hidden layer
\cite{Dwivedi} was proposed for 1D linear advection and diffusion problems.
To decrease the
order of second-order differential operator and the neural network regularity requirement,
the variational PINN was introduced by \cite{Kharazmi2019}. Similarly, inspired by
the classical mixed finite element method and least-squares finite element method,
the deep mixed residual method \cite{Lyu} and first-order deep least-squares
methods \cite{Cai} have been introduced for the general high-order
differential PDEs and second-order elliptic PDEs, respectively.
The asymptotic preserving first-order deep neural network for the
anisotropic elliptic problems was studied in \cite{Li}.
However, the diffusion equations with the general discontinuous
anisotropic tensor are not considered thoroughly in the literature.
The anisotropic heat flux across the discontinuous interface is generally discontinuous and
so do the solution gradient. The derivatives of these variables near the
interface in the loss function will make the training
process more difficult.
Indeed, as pointed in a recent
comprehensive review paper \cite{Karniadakis},
solving PDEs with base PINN indeed may fail in some low-dimensional
cases, such as the diffusion equation with non-smooth conductivity.

Therefore, in this paper,
we introduce a deep neural network solver based on
a new weighted first-order formulation of elliptic anisotropic diffusion problems.
More specifically, the original second-order anisotropic
diffusion equation is rewritten into a
first-order formulation by
introducing two auxiliary variables.
These auxiliary variables
and solution gradient may be discontinuous across the
discontinuity interface. The derivative of these variables is not
well-defined near the interface pointwisely. A weighted function
with the decaying property near the
interface is used to multiplied with the first-order formulation,
thus the training process will not reinforce
a smooth transition nearby the interface and allow to inherit the
true jumps of the heat flux. Therefore,
the accuracy of the network approximation for the solution is much improved.
This is because the anisotropic heat flux is approximated properly
in each smooth domain and correctly handled across
the discontinuity interfaces.

The remainder of this paper is organized as follows. In Section 2, we will give the
problem setting, a weighted first-order formulation. In Section 3, the deep neural
network solver is introduced.
In Section 4, we present some numerical results to
verify the accuracy improvement.
Finally, some conclusions are given in Section 5.
\section{A weighted first-order formulation}
\subsection{The problem setting}
 We consider the following anisotropic diffusion equation in $\Om\subset\R^2$:
 \be
 -\rtri\cdot(\Lam \rtri u)&=&f \text{\ in\ }\Om,\label{equation1} \\
 u&=&g \text{\ on\ }\pa \Om,\label{eq1-4}
 \ee
where $\Lam(x)$ is a $2\times 2$ diffusion tensor, which is a symmetric and
positive definite (SPD) function,
$g(x)$ is the prescribed boundary condition, and the function
$f(x)$ is the so-called source term. If we have $f\geq 0$,
$g\geq0$, then the non-negativity of $u$ is assured by the maximum principle
in the continuum setting.

\subsection{A weighted first-order formulation}
In this section, we will give an equivalent first-order formulation of (\ref{equation1})
-(\ref{eq1-4}).

Let the diffusion tensor $\Lam(x)$ be
\ben
         \Lam=
          \begin{pmatrix}
        a_{11}(x)&a_{12}(x)\\
        a_{12}(x)& a_{22}(x)
          \end{pmatrix}.
\een
In the sequel, we will omit the space-dependence on $x$ of each component $a_{ij}$ of $\Lam$.
Based on the elementary linear algebra, the eigenvalues of $\Lam$ are
\ben
 \lambda_1&=&\frac{a_{11}+a_{22}+\sqrt{(a_{11}-a_{22})^2+4a_{12}^2}}{2}, \\
 \lambda_2&=&\frac{a_{11}+a_{22}-\sqrt{(a_{11}-a_{22})^2+4a_{12}^2}}{2}.
 \een
 With the two eigenvalues, we have two orthogonal eigenvectors as
   \ben
  \hat{q}_1&=&(-a_{12},\lambda_2-a_{22})^T, \\
  \hat{q}_2&=&(\lambda_1-a_{11},-a_{12})^T.
 \een
 Then the unit eigenvectors are given as
  \ben
  q_1&=&\frac{\hat{q}_1}{||\hat{q}_1||}, \\
  q_2&=&\frac{\hat{q}_2}{||\hat{q}_2||}.
 \een
 With these two unit eigenvectors, we can decompose the heat flux $q=\Lambda  \rtri u$ as
 \ben
  q = (\Lambda  \rtri u, q_1)q_1+ (\Lambda  \rtri u, q_2)q_2.
 \een
 Moreover, we define two auxiliary variables as
 \ben
  \tau&=&(\Lambda  \rtri u, q_1), \\
  \phi&=&(\Lambda  \rtri u, q_2),
 \een
 then we have a first-order formulation equivalent to (\ref{equation1})
-(\ref{eq1-4})
  \be
   \left\{
         \begin{array}{ll}
        -\rtri\cdot(\tau q_1+\phi q_2) = f, &\text{\ in\ }\Om,\\
        \tau=(\Lambda  \rtri u, q_1), &\text{\ in\ }\Om,\\
        \phi=(\Lambda  \rtri u, q_2),&\text{\ in\ }\Om, \\
         u=g, &\text{\ on\ }\pa \Om.
         \end{array}
          \right.
   \ee

   Given that the diffusion tensor $\Lam$ may have the discontinuity
   across the interface, the exact solution gradient $\rtri u$ and two
   variables ($\tau$, $\phi$) may be discontinuous.
   The point-wise definition in PINN may be ill-posed in
   the points nearby the interface. Namely, more points
   sampled in this region, will cause the derivative of those discontinuous
   variables become more larger, so do the loss function. Thus the larger
   dispersion of these variables will appear for making the loss function
   decreased. The resulted neural network approximation will deviate the
   true solution profile. To relieve this deviation, we define the distance-weighted
    function $\zeta(x)$ in the whole domain. The function $\zeta(x)$ will
    approach to zero when $x$ is approaching to the interface, and be
    $O(1)$ in the other domains.
     \begin{figure}[!htbp]
\centerline{\includegraphics[width=4cm,
height=4cm]{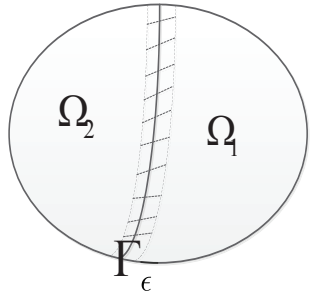}}
\caption{The jumping interface transition region $\Gamma_{\epsilon}$.}
\label{fig1}
\end{figure}
    For convenience, we define the transition
    domain for each interface $\Gamma$ as
    $$\Gamma_\epsilon=\{x\in\Omega, \text{dist}(x, \Gamma)\leq\epsilon\},$$
    then
     \be
   \zeta(x)= \left\{
         \begin{array}{ll}
        \frac{\text{dist}(x,\Gamma^i)}{\text{dist}(x,\Gamma^i)+\epsilon}, &\text{\ in\ }\Gamma^i_\epsilon,\\
         \frac{1}{2}, &\text{\ in\ }\Om\backslash\cup_i\Gamma^i_\epsilon,
         \end{array}
          \right.
   \ee
   where $\epsilon$ is the pre-defined parameter of transition width.
   For example, the domain occupies one diffusion tensor discontinuity interface $\Gamma$
   as in Fig. \ref{fig1}, then the transition domain $\Gamma_\epsilon$ is the shadow
   domain, and $\zeta(x)=\frac{\text{dist}(x,\Gamma)}{\text{dist}(x,\Gamma)+\epsilon}$ in
   this region, equal to $\frac{1}{2}$ otherwise.

With the help of function $\zeta(x)$, we have a new weighted first-order formulation as
 \be\label{weightsystem}
   \left\{
         \begin{array}{ll}
        \zeta(f+\rtri\cdot(\tau q_1+\phi q_2)) = 0, &\text{\ in\ }\Om,\\
        \zeta(\tau-(\Lambda  \rtri u, q_1))=0, &\text{\ in\ }\Om,\\
        \zeta(\phi-(\Lambda  \rtri u, q_2))=0,&\text{\ in\ }\Om, \\
         u=g, &\text{\ on\ }\pa \Om.
         \end{array}
          \right.
   \ee
\begin{remark}
 For 3D anisotropic case, three orthogonal eigenvectors ($q_1$, $q_2$, $q_3$) are
 also available and can be used to define three auxiliary variables ($\tau$, $\phi$, $\psi$).
 The similar distance-weighted function $\zeta$ should be applied.
 In 2D/3D isotropic case, we just choose $q_i=e_i$, which is the
 unit elementary basis in cartesian coordinates, and the width of transition
 domain can be chosen as the machine precision.
\end{remark}
\section{Deep neural network solver}
\subsection{Deep neural network architecture}
\begin{figure}[!htbp]
\centerline{\includegraphics[width=10cm,
height=5cm]{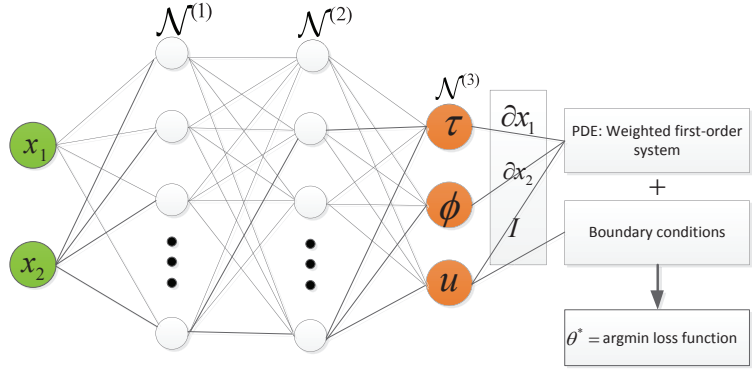}}
\caption{The network structure ($L=3$), the output is then evaluated with
automatic differentiation to get the PDE residual and boundary condition residual.}
\label{fig2}
\end{figure}

We will briefly introduce the deep neural network structure used in this paper.
Generally, a deep neural network defines a nonlinear function
$$\mathcal{N}: \mathbf{x}\in\mathbb{R}^d\rightarrow \mathbf{y}\in\mathbb{R}^c,$$
where $d$ and $c$ are the dimensions of the input and output, respectively.
$\mathcal{N}(\mathbf{x})$ is a composite function of some different layers
of functions as
$$\mathcal{N}(\mathbf{x})=\mathcal{N}^{(L)}\circ\cdots\mathcal{N}^{(2)}\circ\mathcal{N}^{(1)}(\mathbf{x}),$$
where $L$ is the number of layers or the so-called depth of the neural network
(see Fig. \ref{fig2}).

\begin{figure}[!htbp]
\centerline{\includegraphics[width=6cm,
height=4.5cm]{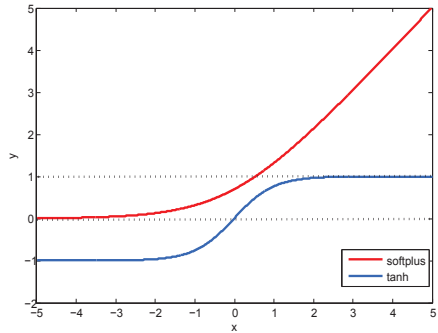}}
\caption{The activation function and softplus function.}
\label{fig3}
\end{figure}

For the $k$th layer $\mathcal{N}^{(k)}$, we use full connected hidden layers, that is,
$$\mathcal{N}^{(k)}(\mathbf{x}^{k-1})=\sigma(W^k\mathbf{x}^{k-1}+b^k),\quad 1\leq k< L,$$
$$\mathcal{N}^{(L)}(\mathbf{x}^{L-1})= W^L\mathbf{x}^{L-1}+b^L, \quad k= L,$$
where $n_k$ is the width of $k$th layer,
the matrix $W^k\in \mathbb{R}^{n_k\times n_{k-1}}$ and vector $b^k\in \mathbb{R}^{n_k}$ are
called the weights and bias, respectively. All the weights and bias
are denoted as $\theta=\{(W^k,b^k)\}_{k=1,\ldots,L}$, which are the
parameters to be trained. The function $\sigma$ is generally nonlinear and called the
activation function. Several classical activation functions are ReLU function,
sigmoid function, hyperbolic tangent function, etc. In this paper, the
hyperbolic tangent function is used as the activation function (see Fig. \ref{fig3})
$$\sigma(\mathbf{x}) = \text{tanh}(\mathbf{x}),$$
where it is applied component-wisely.  With this activation function and multiple
hidden layers, the deep neural network can approximate a large
class of complex nonlinear function rather than the linear ones.

\subsection{The discrete least-squares functional}
In this part, the least-squares formulation for the first-order weighted
formulation (\ref{weightsystem}) is given, then the discrete version based on the
deep neural network approximation is introduced.

The first-order weighted least-squares formulation is to find $(u,\tau,\phi)\in (H^1(\Omega))^3$
such that
\be
\Psi(u,\tau,\phi;\mathbf{f})= \min_{(\eta,\nu,\chi)\in (H^1(\Omega))^3}
                               \Psi(\eta,\nu,\chi;\mathbf{f}),
\ee
where $\mathbf{f}=(f,g,\zeta)$ and
\be
\Psi(\eta,\nu,\chi;\mathbf{f})&=& \|\zeta(f+\rtri\cdot(\nu q_1+\chi q_2))\|_{0,\Om}^2
                               +\|\zeta(\nu-(\Lambda  \rtri \eta, q_1))\|_{0,\Om}^2\nn\\
                               &&+\| \zeta(\chi-(\Lambda  \rtri \eta, q_2))\|_{0,\Om}^2
                               +\|\eta-g\|_{1/2,\partial\Omega}^2.
\ee
 If three unknown functions $u$, $\tau$, $\phi$ are approximated by one deep
 neural network and its three outputs are
 denoted by $\hat{u}(x,\theta)$, $\hat{\tau}(x,\theta)$, $\hat{\phi}(x,\theta)$ (see Fig. \ref{fig2}), then the discrete formulation based on all sampling points reads
 \be
\hat{\Psi}(\hat{u},\hat{\tau},\hat{\phi};\mathbf{f})(\theta)= \min_{\tilde{\theta}\in\mathbb{R}^N}
                               \hat{\Psi}(\hat{\eta},\hat{\nu},\hat{\chi};\mathbf{f})(\tilde{\theta}),
\ee
where the discrete functional reads
\be
\hat{\Psi}(\hat{\eta},\hat{\nu},\hat{\chi};\mathbf{f})(\tilde{\theta})&=&
\frac{1}{N_f}\sum_{i=1}^{N_f}(\zeta(\mathbf{x}_i)(f(\mathbf{x}_i)+\rtri\cdot
(\hat{\nu}(\mathbf{x}_i,\tilde{\theta})
 q_1(\mathbf{x}_i)+\hat{\chi}(\mathbf{x}_i,\tilde{\theta}) q_2(\mathbf{x}_i))))^2\nn\\
 &&+\frac{1}{N_f}\sum_{i=1}^{N_f}(\zeta(\mathbf{x}_i)(\hat{\nu}(\mathbf{x}_i,\tilde{\theta})-(\Lambda  \rtri \hat{\eta}(\mathbf{x}_i,\tilde{\theta}), q_1(\mathbf{x}_i))))^2\nn\\
 &&+ \frac{1}{N_f}\sum_{i=1}^{N_f}(\zeta(\mathbf{x}_i)(\hat{\chi}(\mathbf{x}_i,\tilde{\theta})-(\Lambda  \rtri \hat{\eta}(\mathbf{x}_i,\tilde{\theta}), q_2(\mathbf{x}_i))))^2\nn\\
 &&+\frac{\omega_D}{N_D}\sum_{i=1}^{N_D}(\hat{\eta}(\mathbf{x}_i,\tilde{\theta})-g(\mathbf{x}_i))^2,
 \ee
 where $N_f$, $N_D$ are number of collocation points in $\Omega$ and on $\partial\Omega$,
 the total number of network parameters is $N=\sum_{i=1}^L n_i(n_{i-1}+1)$ and
 $\omega_D$ represents the boundary-weight to penalize the neural
 network approximations satisfying the boundary conditions.
 \begin{remark}
As pointed in deep mixed residual method \cite{Lyu} and deep least-squares method
\cite{Cai}, the three unknown functions can also
employ more than one neural network (multiple different neural networks)
to approximate them. This may offer another flexibility
and accuracy improvement for solving this weighted first-order system.
\end{remark}
\begin{remark}
In some cases, if the function $\mathbf{f}\ge 0$, the analytical solution $u$ is non-negative
by the maximum principle. To enforce this constraints, two alternative ways are possible.
One is the penalty way similar to the boundary condition enforcement reads
$$\hat{\Psi}_{p_1}(\hat{\eta},\hat{\nu},\hat{\chi};\mathbf{f})(\tilde{\theta})=
\hat{\Psi}(\hat{\eta},\hat{\nu},\hat{\chi};\mathbf{f})(\tilde{\theta})
+\frac{\omega_{p_1}}{N_f}\sum_{i=1}^{N_f}\max(-u(\mathbf{x}_i,\tilde{\theta}),0)^2,$$
where $\omega_{p_1}$ is the penalty parameter.
Another way is to apply the softplus function $\ln(1+e^x)$ (see Fig. \ref{fig3}) to
the network output $\hat{\eta}$ as $\bar{\eta}=\ln(1+e^{\hat{\eta}})$, then
the discrete functional is replaced as
$$\hat{\Psi}_{p_2}(\bar{\eta},\hat{\nu},\hat{\chi};\mathbf{f})(\tilde{\theta})=
\hat{\Psi}(\bar{\eta},\hat{\nu},\hat{\chi};\mathbf{f})(\tilde{\theta}),$$
and then the final neural network approximation becomes $(\bar{u},\tau,\phi)$ with
$\bar{u}\ge 0$ almost everywhere.
\end{remark}
\begin{remark}
For completeness and comparison, we also give the base PINN discrete
formulation for the anisotropic diffusion problems as follows
\be
\hat{\Psi}_{PINN}(\hat{u};\mathbf{f_1})(\theta)= \min_{\tilde{\theta}\in\mathbb{R}^N}
                               \hat{\Psi}_{PINN}(\hat{\eta};\mathbf{f_1})(\tilde{\theta}),
\ee
where the discrete functional reads
\be
\hat{\Psi}_{PINN}(\hat{\eta};\mathbf{f_1})(\tilde{\theta})&=&
\frac{1}{N_f}\sum_{i=1}^{N_f}((f(\mathbf{x}_i)+\rtri\cdot
(\Lambda(\mathbf{x}_i)\rtri \hat{\eta}(\mathbf{x}_i,\tilde{\theta}))))^2\nn\\
 &&+\frac{\omega_D}{N_D}\sum_{i=1}^{N_D}(\hat{\eta}(\mathbf{x}_i,\tilde{\theta})-g(\mathbf{x}_i))^2,
 \ee
 where $\mathbf{f_1}=(f,g)$, $N_f$, $N_D$ are number of collocation points
 in $\Omega$ and on $\partial\Omega$, and
 $\omega_D$ represents the boundary-weight to penalize the neural
 network approximations satisfying the boundary conditions.
\end{remark}

\section{Numerical results}

In this section,  we test several diffusion problems with
continuous or discontinuous anisotropic
diffusion tensors. For simplicity, our method is denoted as weighted FO-PINN.
For both base PINN and our weighted FO-PINN,
we use the same neural network structure.
More specifically, we use a five-layer neural network with 16 neurons
to approximate the solution
$u$ and two auxiliary variables simultaneously in our method, and a five-layer
neural network with 16 neurons just for the solution approximation in base PINN method.
The activation function for all
hidden layers is the tanh function and identity for the last output layer.

If there exists the analytical solution,
we use the following relative discrete $L_2
$-norm to evaluate the accuracy of learned solution
\ben
\epsilon_2^u&=&\frac{\|u-u_e\|_2}{\|u_e\|_2},\\
\een
where $u_e$ is the problem-dependent analytical solution.

To implement and train the neural network, an integrated Python software
deepXDE \cite{Lulu} based TensorFlow is used to implement the base PINN and
our method. To train the neural network parameters $\theta$ efficiently, the Adam optimizer version
of gradient descent is used with a pre-defined number of iterations, then the
L-BFGS-B optimizer is used for achieving the final convergence. The learning
rate is initialized as $1e-3$ and will decrease gradually
as the iteration steps growing. The collocation points are sampled randomly from
the interior domain and domain boundary with the specific numbers $N_f$ and $N_D$.

\subsection{The diffusion problem with smooth anisotropic coefficient}
\label{smooth}
We first test the accuracy for diffusion problem
with smooth anisotropic coefficients. Consider the problem (\ref{equation1}) with Dirichlet
boundary condition in the unit square $\Omega =(0,  1)^2$.
The exact analytical solution and the diffusion coefficient coming from
the FVCA6 benchmark \cite{Eymard} are
$u=1+\sin (2\pi x)\sin (2\pi y)$ and \[\Lambda =\left[
\begin{matrix}
   1 & {}  \\
    {} & 10^3  \\
\end{matrix}\right].
\]

The source $f$ and the Dirichlet boundary data $g$ are set
accordingly to the exact solution.

In this case, we use $N_f=10000$ collocation points, $N_D=1000$ boundary
points and boundary-weight $\omega_D=100$, the Adam optimizer
is used until 10000 iterations and then L-BFGS-B
optimizer is used to reach convergence.
\begin{figure}[!htbp]
\begin{minipage}{0.48\linewidth}
\centerline{\includegraphics[width=5.5cm,height=5cm]{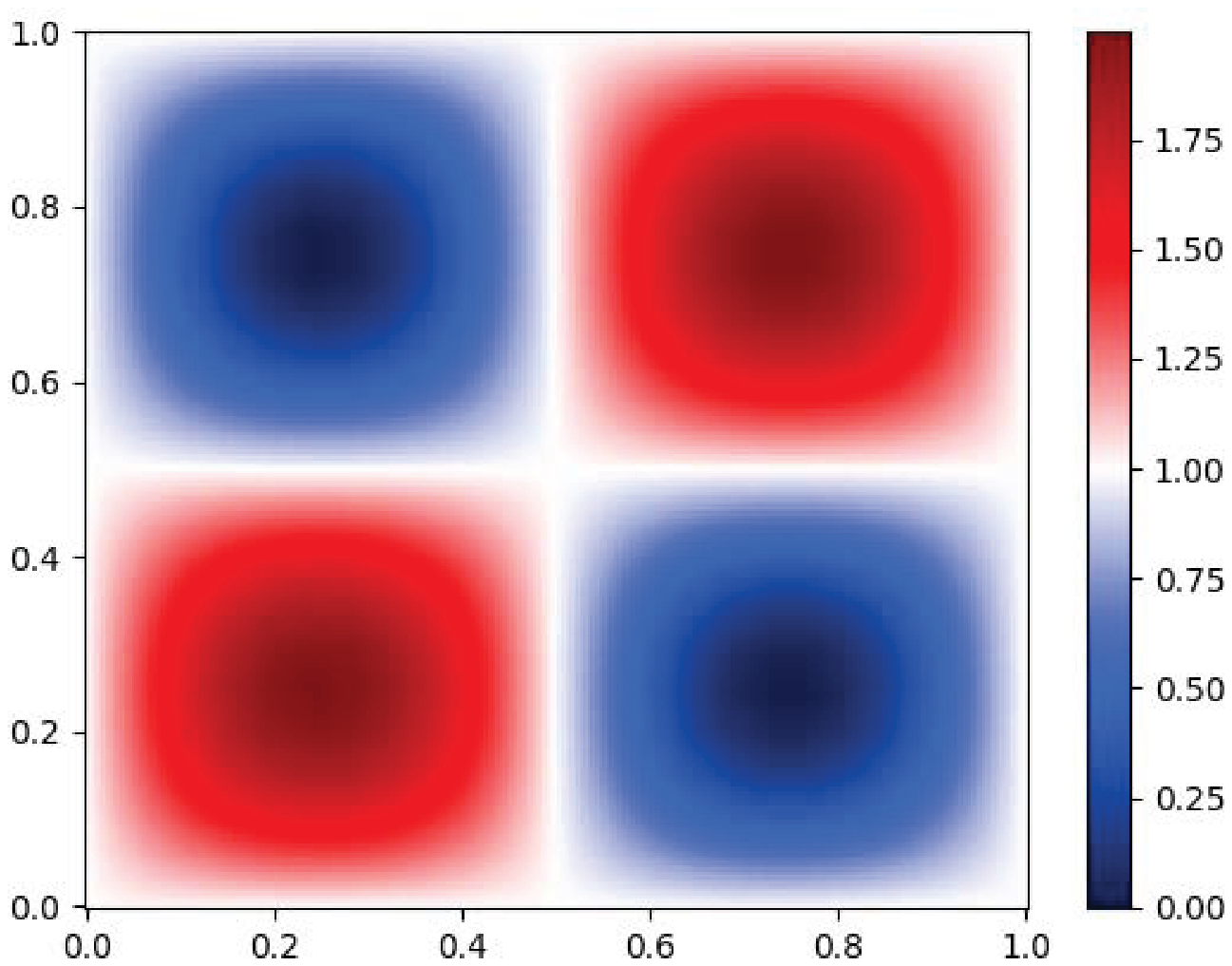}}
\end{minipage}
\hfill
\begin{minipage}{0.48\linewidth}
\centerline{\includegraphics[width=5.5cm,height=5cm]{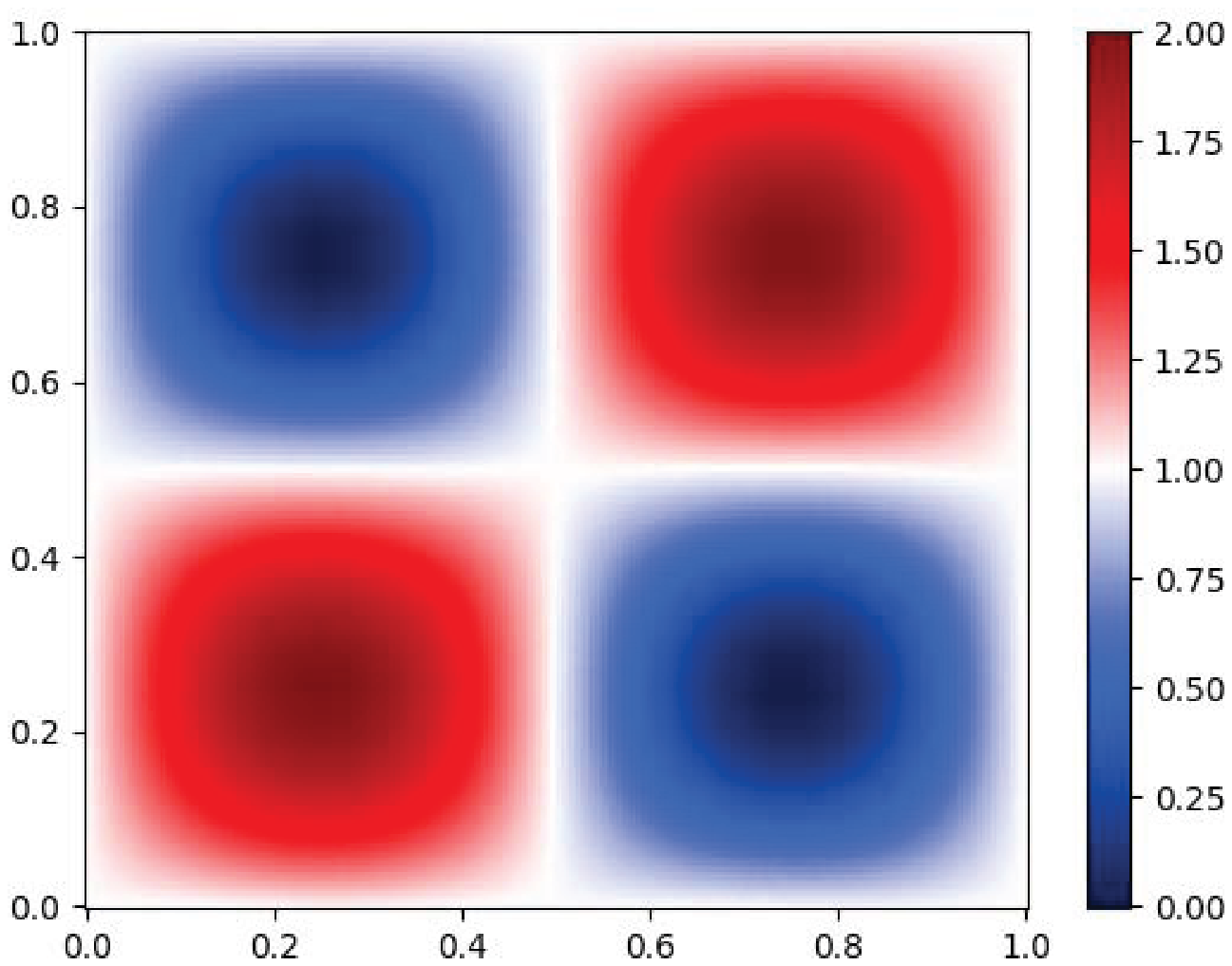}}
\end{minipage}
\caption{The solution profiles of base PINN (left) and
weighted FO-PINN (right). }
\label{fig4}
\end{figure}

The numerical results for base PINN and our weighted FO-PINN
are shown in Fig. \ref{fig4}. The discrete relative $L_2
$-norm for base PINN and weighted FO-PINN are
2.93e-3 and 4.1e-3, respectively.
We can see that both methods are quite accurate for
this smooth problem, and our method can naturally
solve this problem by setting the weight function $\zeta$
equal to $1/2$ in the whole domain.

\subsection{The diffusion problem with discontinuous isotropic coefficient}
Next, we consider the problem (\ref{equation1}-\ref{eq1-4})
with Dirichlet boundary condition in the unit domain
 $\Om=[0,1]^2$. The diffusion tensor $\Lam(x)$ is
 a isotropic but discontinuous tensor and given by
 \ben
        \Lam=\left\{
         \begin{array}{c}
        5,\quad x\leq0.5,\\
        1,\quad x>0.5.\\
         \end{array}
         \right.
 \een
 The analytical solution is chosen to be
 \ben
         u=\left\{
         \begin{array}{c}
        (x^2+10)(y-y^2),\quad x\leq0.5,\\
        (5x^2+9)(y-y^2),\quad x>0.5.\\
         \end{array}
         \right. \nn
 \een
The source $f$ and the Dirichlet boundary data $g$ are set
accordingly to the exact solution.

In this case, we use $N_f=10000$ collocation points, $N_D=2000$ boundary
points, and the interior weights and boundary-weight are
chosen as the inverse of the initial losses of each
term such that the weighted losses are of the
same order of magnitude $O(1)$, the Adam optimizer
is used until $10000$ iterations and then L-BFGS-B
optimizer is used to reach convergence.

 \begin{figure}[!htbp]
\centerline{\includegraphics[width=6.5cm,
height=5cm]{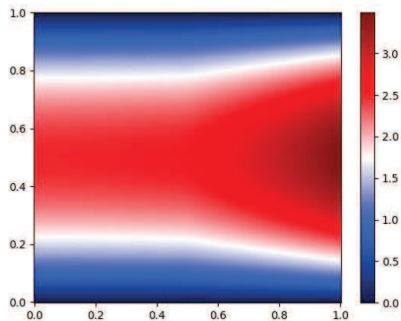}}
\caption{The exact solution profile for isotropic but discontinuous diffusion
tensor.}
\label{fig5-1}
\end{figure}
\begin{figure}[!htbp]
\begin{minipage}{0.48\linewidth}
\centerline{\includegraphics[width=6.5cm,height=5cm]{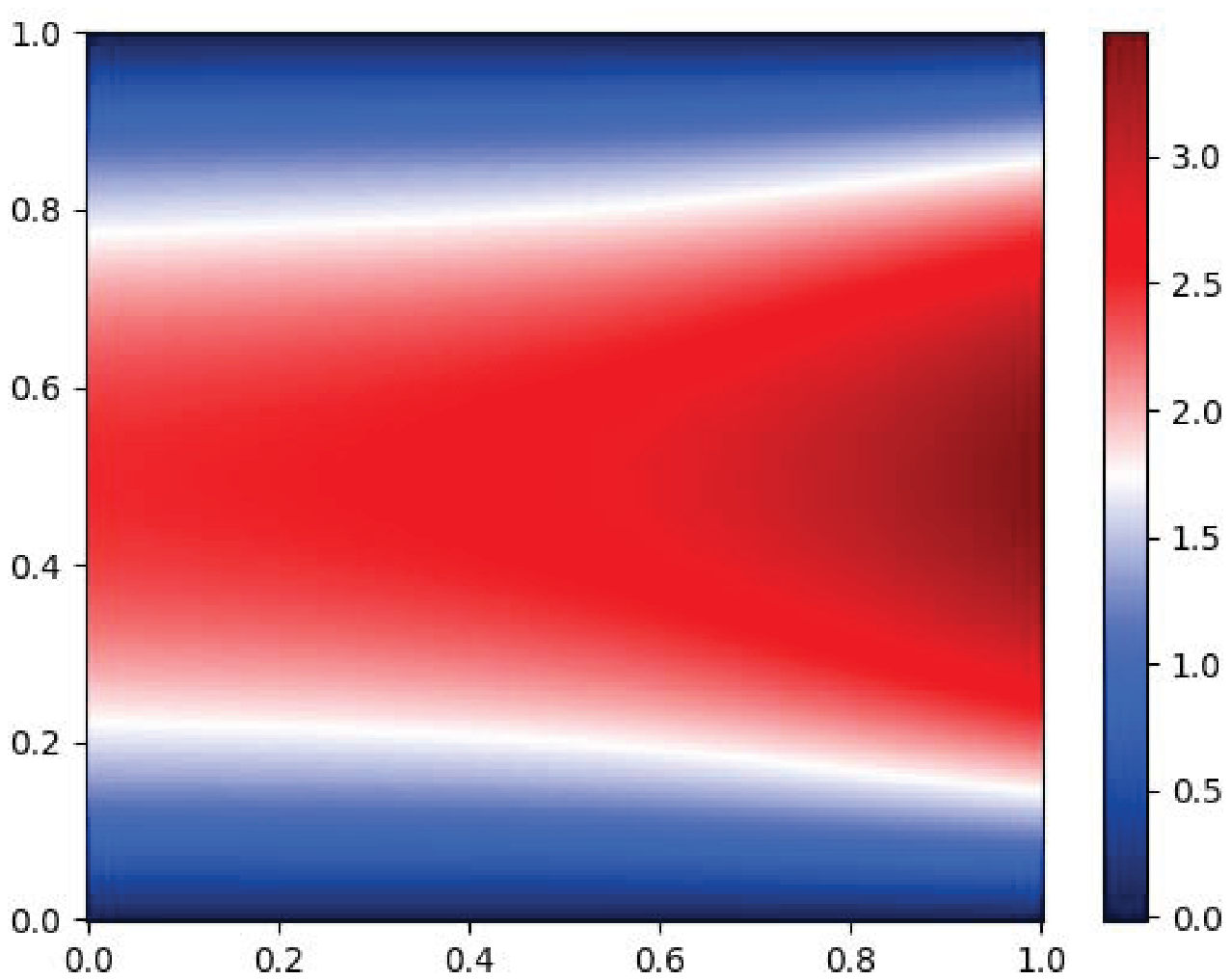}}
\end{minipage}
\hfill
\begin{minipage}{0.48\linewidth}
\centerline{\includegraphics[width=6.5cm,height=5cm]{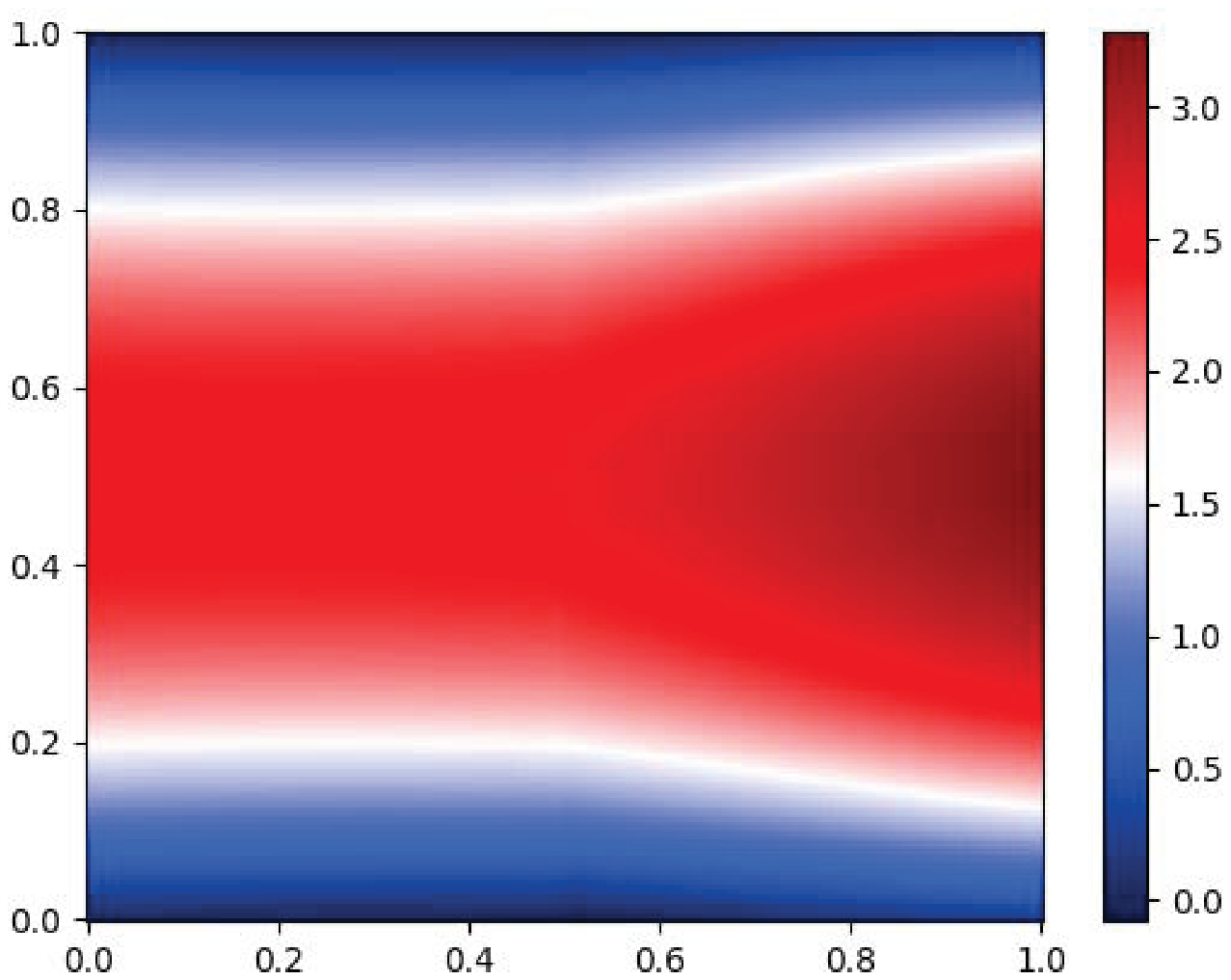}}
\end{minipage}
\caption{The solution profiles of base PINN (left) and
weighted FO-PINN (right). }
\label{fig5-2}
\end{figure}

The exact solution is depicted in Fig. \ref{fig5-1} for comparison.
The numerical results for base PINN and our weighted FO-PINN
are shown in Fig. \ref{fig5-2}. The discrete relative $L_2
$-norm for base PINN and weighted FO-PINN are
$2.64e-2$ and $8.06e-3$, respectively.
We can see that our method is more accurate than
the base PINN method for this discontinuous isotropic case,
and the transition region width $\epsilon$ in the weighted function $\zeta$
is equal to $0.01$. It is shown that the resolution of our method is higher
with the curved contours near the discontinuous interface.

\subsection{The diffusion problem with discontinuous anisotropic coefficient}

In this subsection, we demonstrate the accuracy for diffusion problem with
the discontinuous and mild anisotropic coefficient, where
many traditional cell-centered finite volume
schemes also will undergo the accuracy degeneracy.
Let the exact solution and the diffusion coefficient be as follows:
 \ben
    \Lam=\left\{
         \begin{array}{c}
        \left(\begin{array}{cc}
        1&0\\
        0&1
        \end{array}\right),\quad\ x\leq0.5,\\
       \left(\begin{array}{ccc}
        10&3\\
        3&1
         \end{array}\right),\quad x>0.5,\\
         \end{array}
         \right.
 \een
 and the analytical solution is chosen to be
 \ben
         u(x,y)=\left\{
         \begin{array}{l}
        -2y^2+4xy+6x+2y+1,\quad\quad\qquad x\leq0.5,\\
        -2y^2+1.6xy-0.6x+3.2y+4.3,\quad x>0.5.\\
         \end{array}
          \right.
 \een
The source term and Dirichlet boundary condition are computed by this analytical solution.
\begin{figure}[!htbp]
\begin{minipage}{0.48\linewidth}
\centerline{\includegraphics[width=6.5cm,height=5cm]{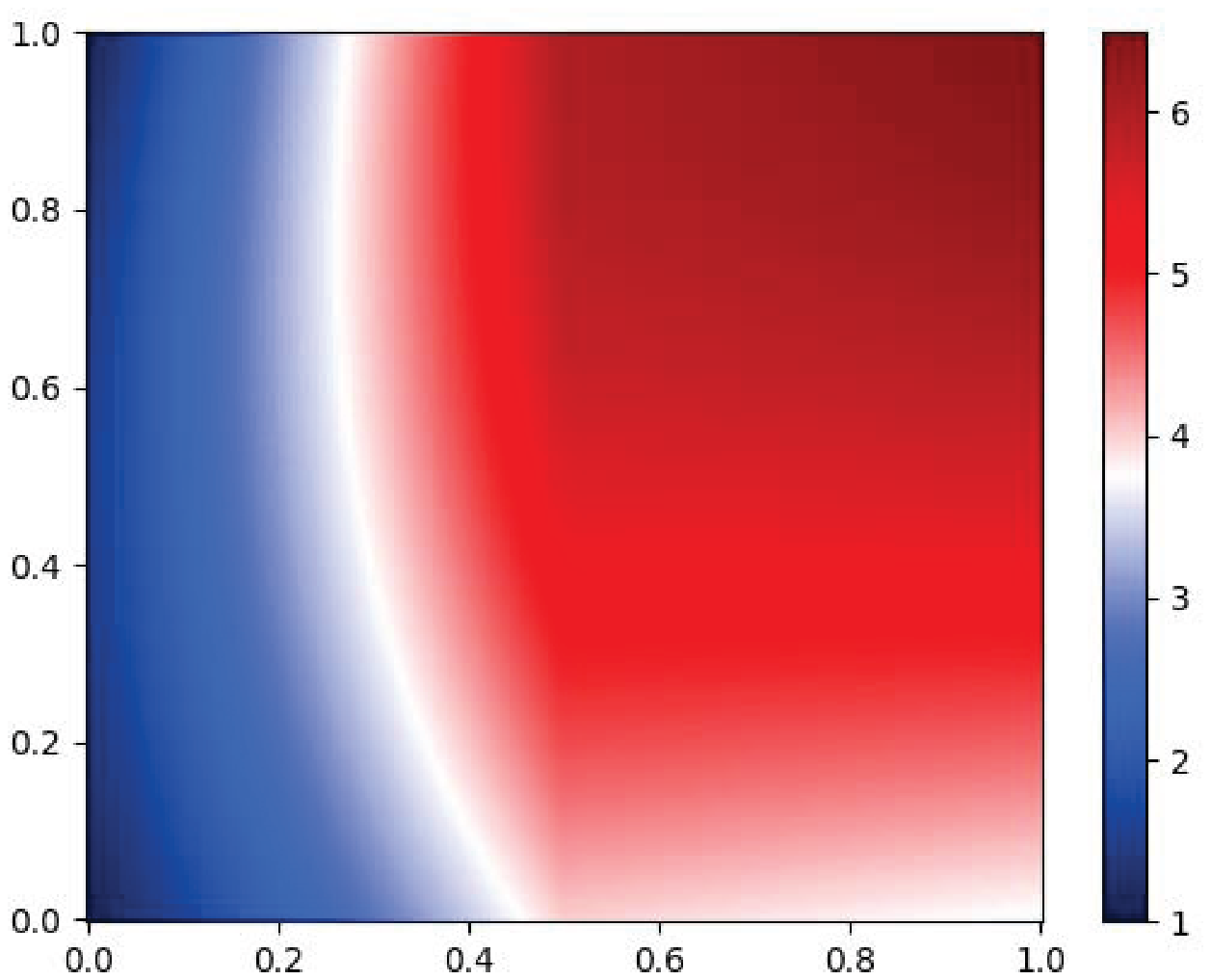}}
\end{minipage}
\hfill
\begin{minipage}{0.48\linewidth}
\centerline{\includegraphics[width=6.5cm,height=5cm]{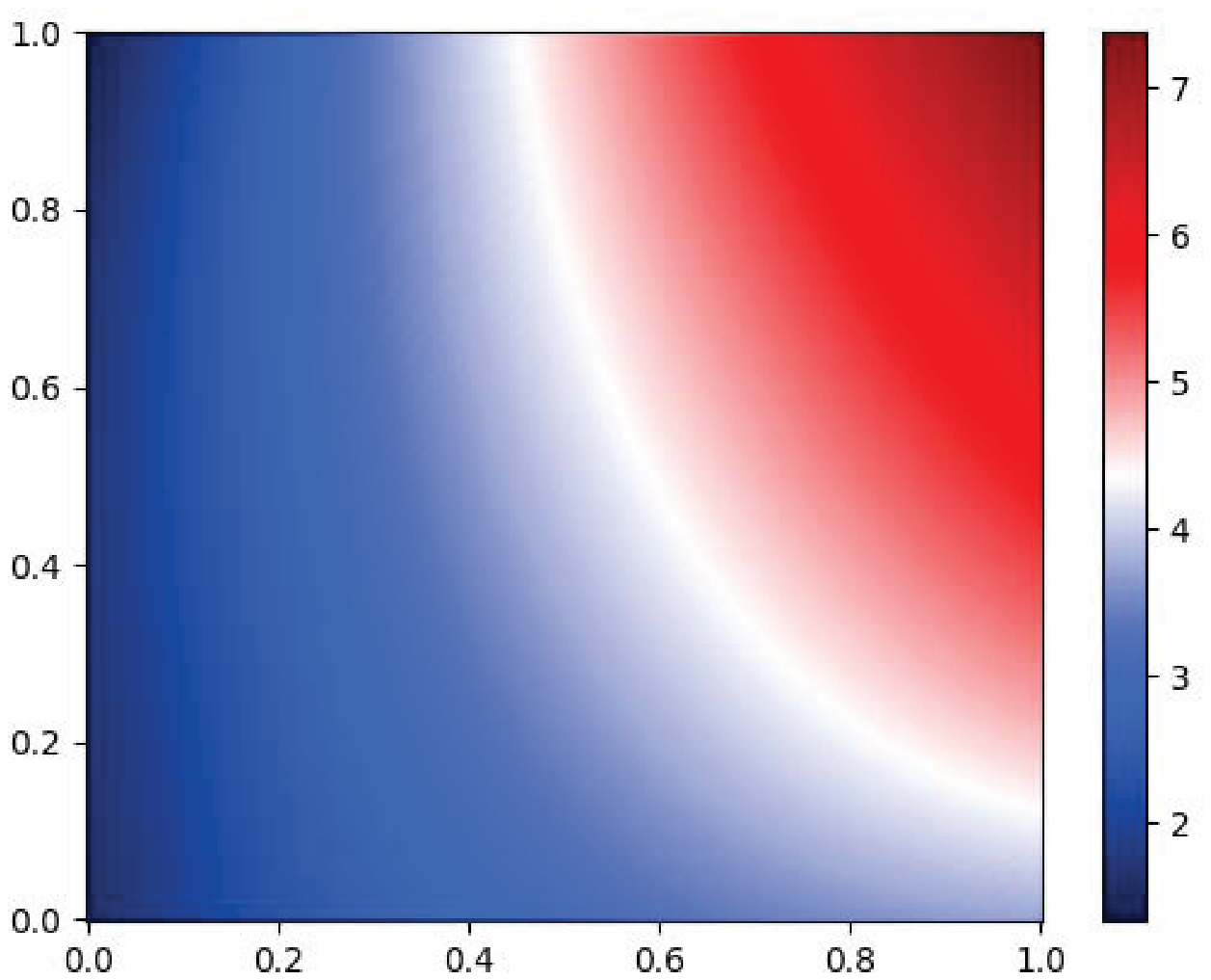}}
\end{minipage}
\caption{The exact solution profile (left) and
the solution profile of base PINN (right). }
\label{fig6-1}
\end{figure}
\begin{figure}[!htbp]
\begin{minipage}{0.48\linewidth}
\centerline{\includegraphics[width=6.5cm,height=5cm]{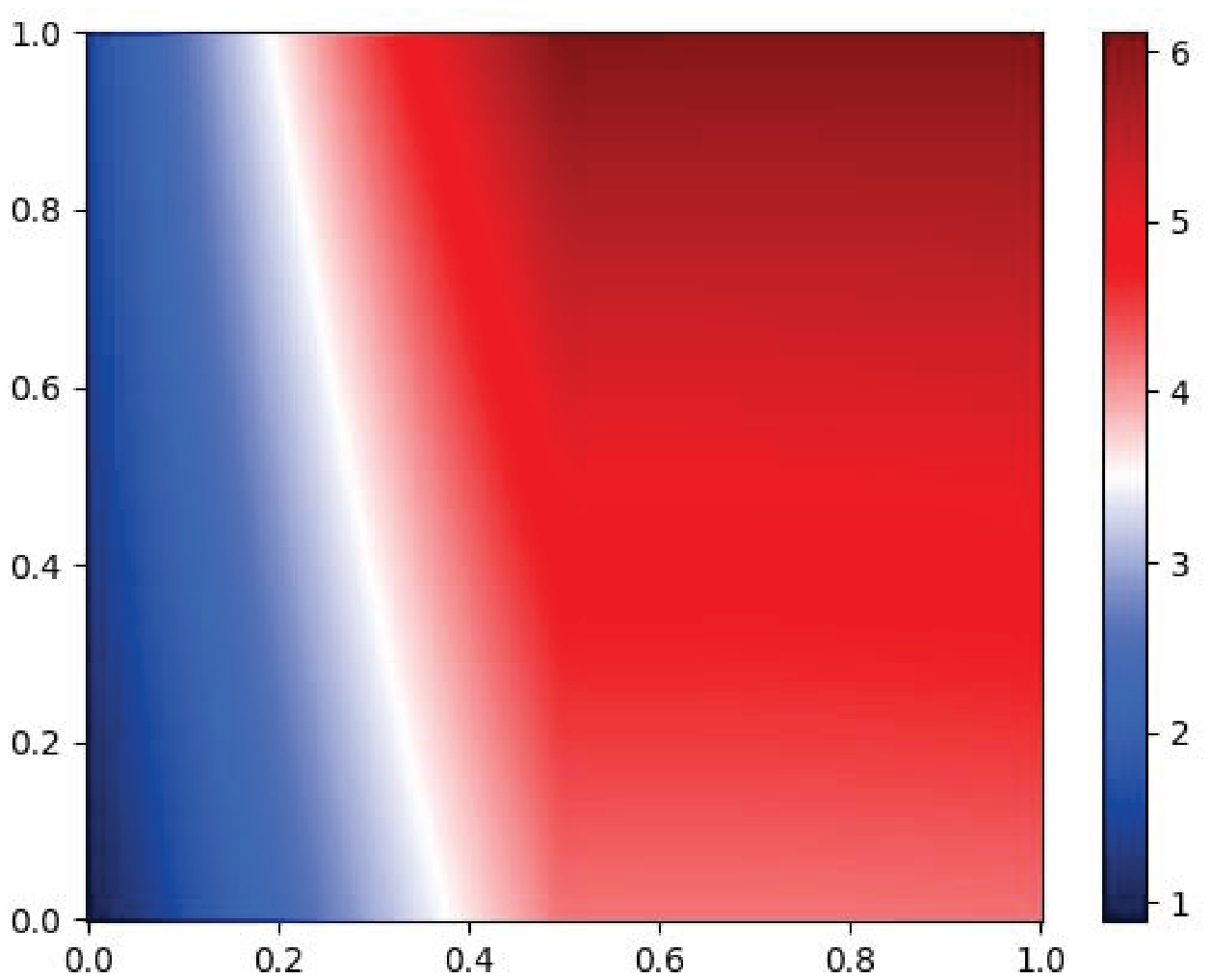}}
\end{minipage}
\hfill
\begin{minipage}{0.48\linewidth}
\centerline{\includegraphics[width=6.5cm,height=5cm]{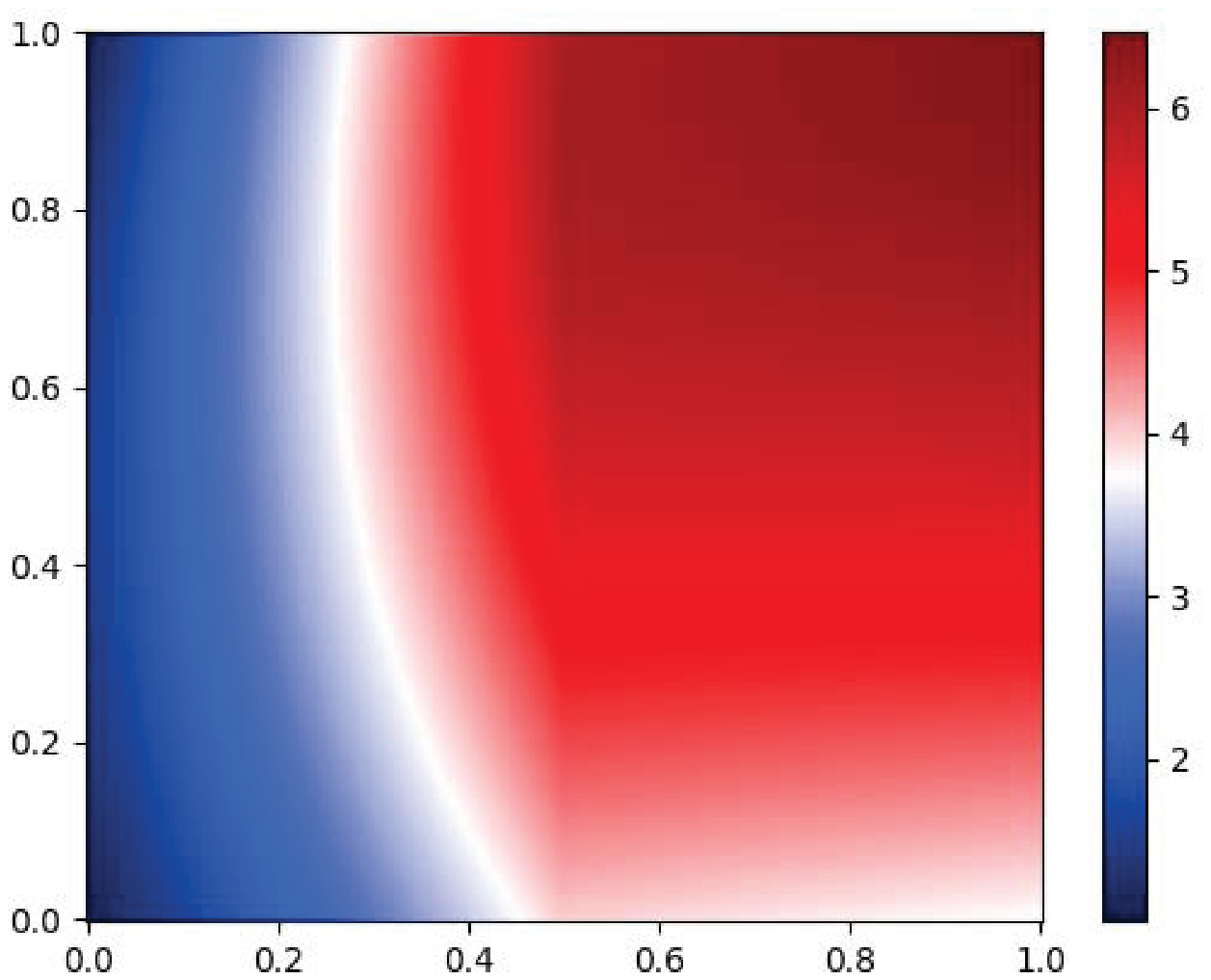}}
\end{minipage}
\caption{The solution profiles of pure FO-PINN (left) and
weighted FO-PINN (right). }
\label{fig6-2}
\end{figure}

In this case, we use the same neural network structure as the last one.
The number of collocation points and
optimization method for training are also
the same as it. The exact solution and the solution of
base PINN method are depicted in Fig. \ref{fig6-1} for comparison.
We can see that the base PINN can not capture the
correct heat flux even in one-side domain.
As shown in Fig. \ref{fig6-2} (left), the pure first-order formulation
is better for capturing the heat flux in one-side domain,
but is less accurate across the jump interface.
By adding the weighted function $\zeta$ with transition
width $\epsilon=0.01$,
The numerical result for the weighted FO-PINN
is shown in Fig. \ref{fig6-2} (right), whose solution profile
is closest to the exact solution.
The discrete relative $L_2
$-norms for base PINN, pure FO-PINN and
weighted FO-PINN are
$1.42e-1$, $7.18e-2$ and $6.31e-3$, respectively.
We can see that our method is much more accurate than
the base PINN method for this discontinuous anisotropic case,
and the weighted function $\zeta$
is crucial for the resolution near the discontinuous interface and
approximation of the heat flux across the interface.

To test the accuracy and efficiency of our method in higher dimension,
we consider a similar 3D case as follows:
 \ben
    \Lam=\left\{
         \begin{array}{c}
        \left(\begin{array}{ccc}
        1&0&0\\
        0&1&0\\
        0&0&1
         \end{array}\right),\quad\ x\leq0.5,\\
       \left(\begin{array}{ccc}
        10&3&0\\
        3&1&0\\
        0&0&1
         \end{array}\right),\quad x>0.5,\\
         \end{array}
         \right.
 \een
 and the analytical solution is chosen to be
 \ben
         u(x,y,z)=\left\{
         \begin{array}{l}
        -2y^2+4xy+6x+2y+z+1,\quad\quad\qquad x\leq0.5,\\
        -2y^2+1.6xy-0.6x+3.2y+z+4.3,\quad x>0.5.\\
         \end{array}
          \right.
 \een
The source term and Dirichlet boundary condition are computed by this analytical solution.
\begin{figure}[!htbp]
\centerline{\includegraphics[width=6.5cm,
height=5cm]{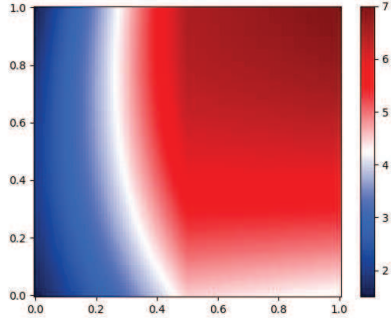}}
\caption{The exact solution slide on $z=0.5$ plane for 3D discontinuous diffusion
tensor.}
\label{fig10-1}
\end{figure}
\begin{figure}[!htbp]
\begin{minipage}{0.48\linewidth}
\centerline{\includegraphics[width=6.5cm,height=5cm]{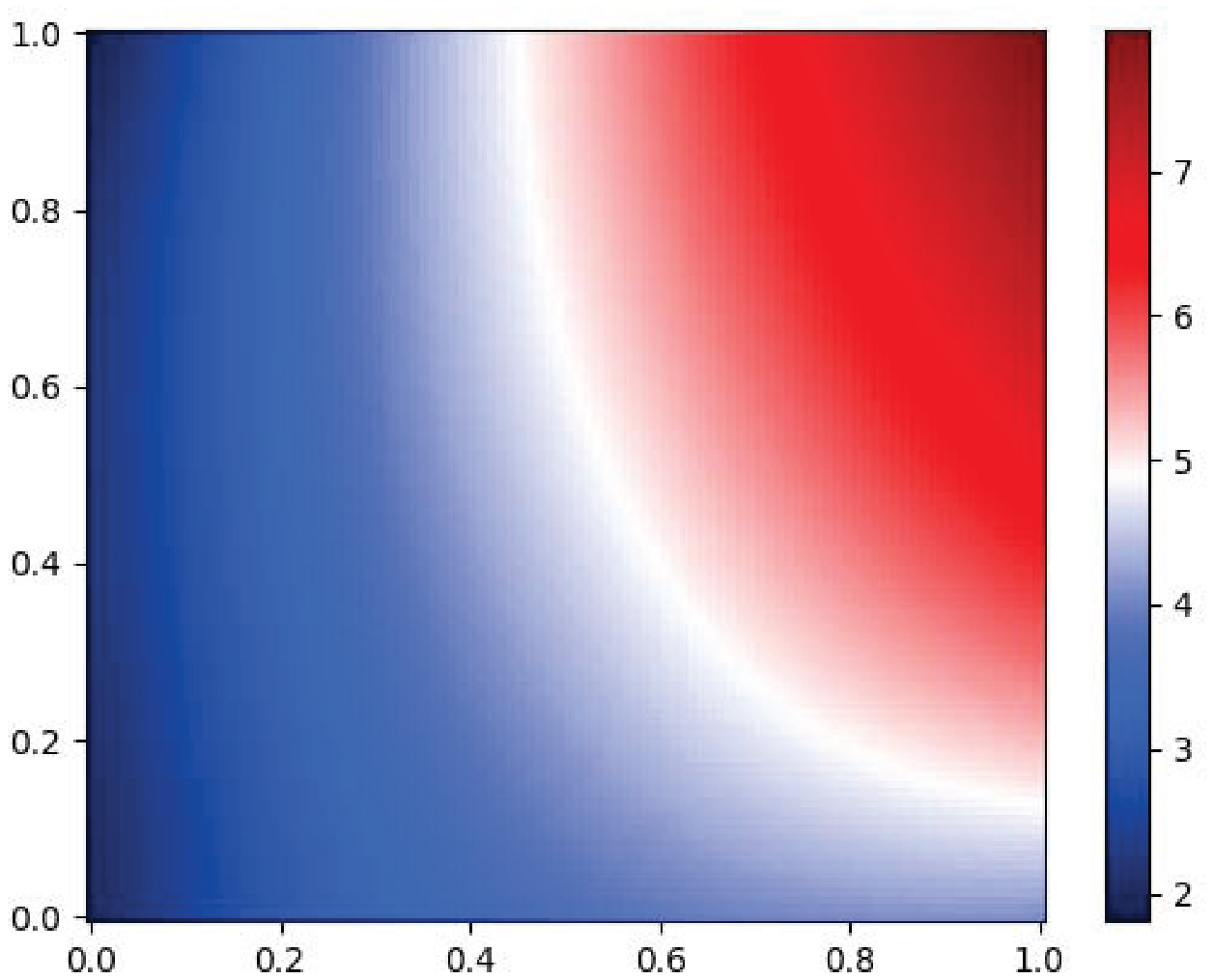}}
\end{minipage}
\hfill
\begin{minipage}{0.48\linewidth}
\centerline{\includegraphics[width=6.5cm,height=5cm]{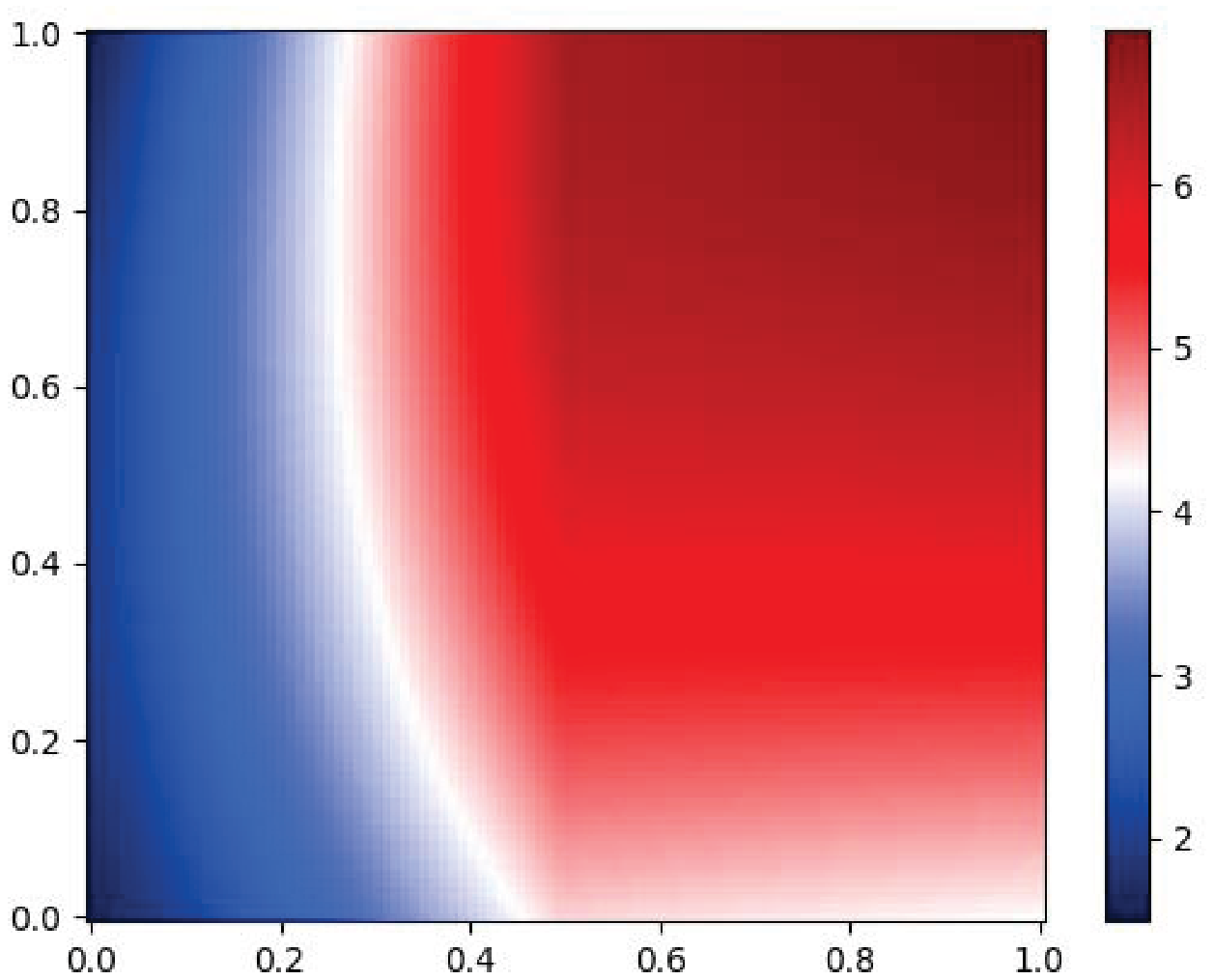}}
\end{minipage}
\caption{The 3D solution slide on $z=0.5$ plane of base PINN (left) and
weighted FO-PINN (right). }
\label{fig10-2}
\end{figure}

In this case, the same neural network structure as 2D case is used, only the
number of outputs is different.
We choose $N_f=20000$ collocation points and $N_D=2000$ boundary
points as twice as 2D counterparts.
The exact solution slide on $z=0.5$ plane is depicted
in Fig. \ref{fig10-1} for comparison.
The numerical results for base PINN and our weighted FO-PINN
are shown in Fig. \ref{fig10-2}. The discrete relative $L_2
$-norm for base PINN and weighted FO-PINN are
$1.19e-1$ and $2.67e-3$, respectively.
We can see that our method is also more accurate than
the base PINN method for this 3D anisotropic case, and the
accuracy of base PINN is slightly worse than 2D, while ours is
more accurate than 2D counterpart. It is shown that
our neural network approximation is less sensitive to the dimensionality.

\subsection{The diffusion problem with strong anisotropic coefficient}

In this subsection, we demonstrate the accuracy
for diffusion problem with a strong anisotropic diffusion coefficient with
more discontinuities. This example is modified from
\cite{Edwards}, such that the gradient of the solution has discontinuity across
the interfaces. The domain $\Omega = [0,1]^2$ is split into four
subdomains $\Omega = \bigcup _{i=1}^4 \Omega_i$, which are given by
\begin{eqnarray*}
  \Omega_1 &=& \{(x,y) \in [0,1]^2 \quad \hbox{such that}  \quad x \leq 0.5, y \leq 0.5  \}\\
  \Omega_2 &=& \{(x,y) \in [0,1]^2 \quad \hbox{such that}  \quad x  > 0.5, y \leq 0.5  \}\\
  \Omega_3 &=& \{(x,y) \in [0,1]^2 \quad \hbox{such that}  \quad x  > 0.5, y  > 0.5  \}\\
  \Omega_4 &=& \{(x,y) \in [0,1]^2 \quad \hbox{such that}  \quad x \leq 0.5, y > 0.5 \}
\end{eqnarray*}

The diffusion coefficient and the exact solution are given by
\[\Lam
=\left[
\begin{matrix}
   a_x^i & 0  \\
   0 & a_y^i
\end{matrix} \right], \ \hbox{for} (x,y) \in \Omega_i, \hbox{ with}\  \begin{tabular}{|c|c|c|c|c|}
                                                                        \hline
                                                                        i    & 1    & 2   & 3    & 4 \\
                                                                        \hline  $a_x^i$    & 10    & 0.1   & 0.01   & 100\\
                                                                        \hline $a_y^i$    & 0.01   & 100 & 10 & 0.1\\                                                                     \hline  $\alpha_i$ & 0.1  & 10  & 100  & 0.01 \\
                                                                        \hline
                                                                        \end{tabular}
\] and \[
 u(x,y) = \alpha_i \sin(2\pi x)\sin(2\pi y).
\]

\begin{figure}[!htbp]
\centerline{\includegraphics[width=6.5cm,
height=5cm]{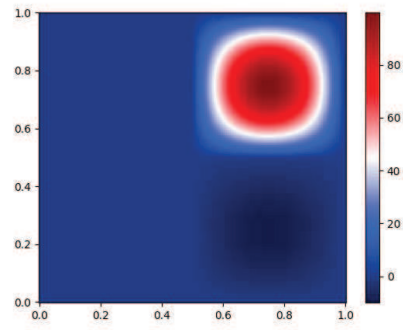}}
\caption{The exact solution profile for strong discontinuous diffusion
tensor.}
\label{fig7-1}
\end{figure}
\begin{figure}[!htbp]
\begin{minipage}{0.48\linewidth}
\centerline{\includegraphics[width=6.5cm,height=5cm]{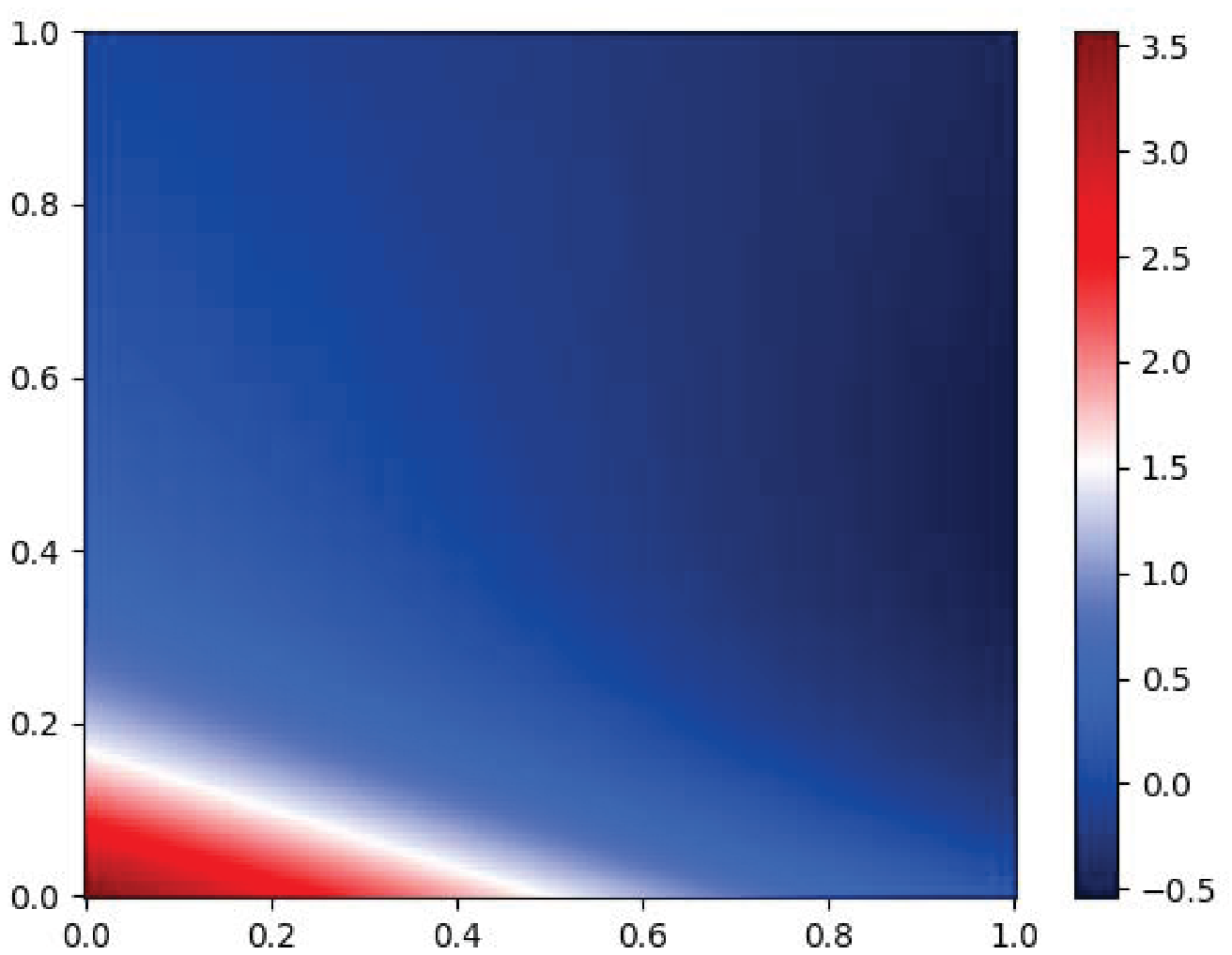}}
\end{minipage}
\hfill
\begin{minipage}{0.48\linewidth}
\centerline{\includegraphics[width=6.5cm,height=5cm]{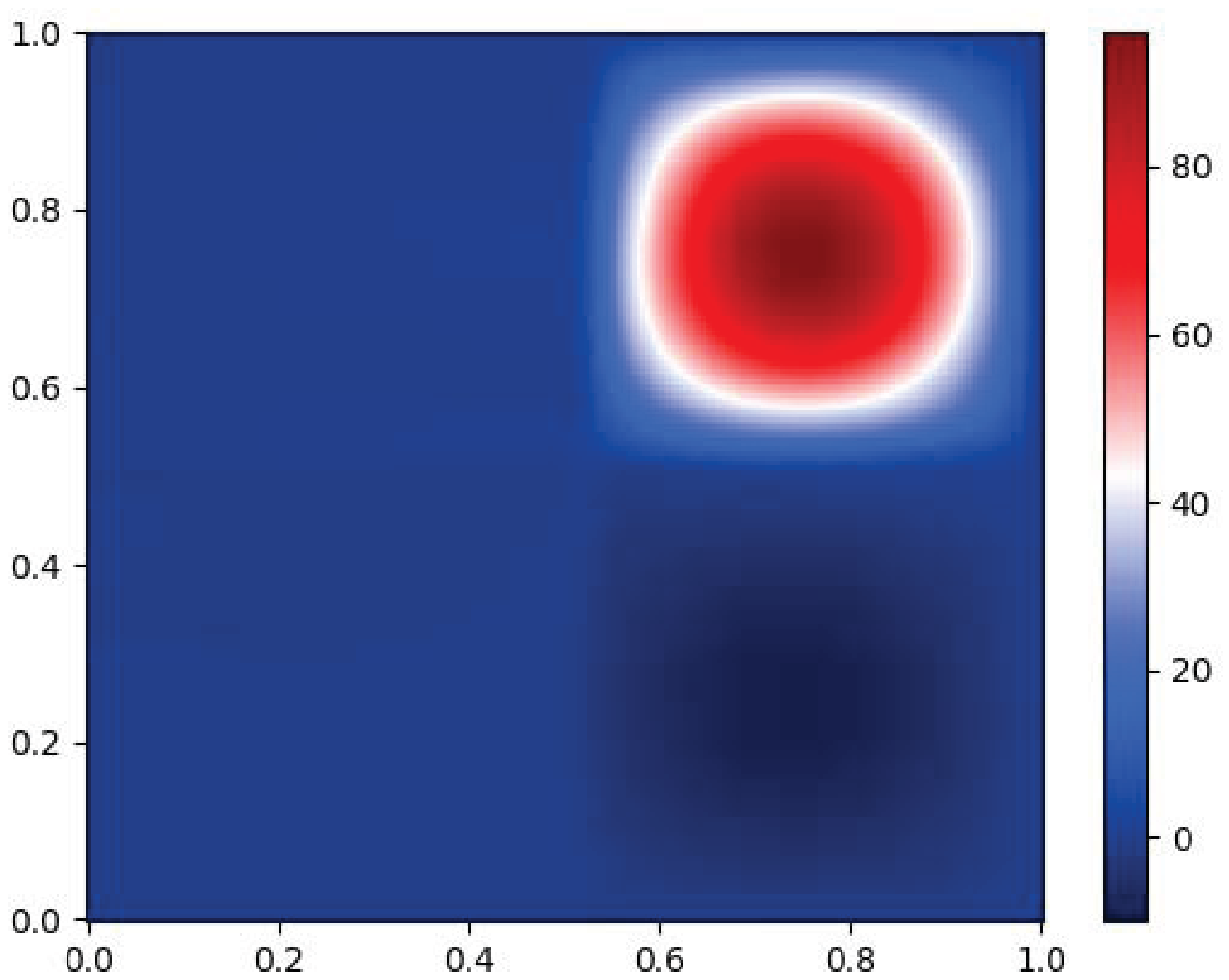}}
\end{minipage}
\caption{The solution profiles of base PINN (left) and
weighted FO-PINN (right). }
\label{fig7-2}
\end{figure}

In this case, we use the same neural network structure as the last one.
The number of collocation points and
optimization method for training are also
the same as it. The exact solution are depicted in
Fig. \ref{fig7-1} for comparison.
As shown in Fig. \ref{fig7-2} (left), the base
PINN is failed to capture the solution
profile, while our weighted FO-PINN
is quite close to the true solution, see
Fig.  \ref{fig7-2} (right).
The discrete relative $L_2
$-norms for base PINN and
weighted FO-PINN are
$1.01e-0$ and $6.90e-2$, respectively.
We can see that our method is much more accurate than
the base PINN method for this strong discontinuous anisotropic case,
and the weighted function $\zeta$
with parameter $\epsilon=0.01$ is used
for the resolution near the discontinuous interface and
approximation of the heat flux across the internal interfaces.

\subsection{Anisotropic diffusion problems with local source term}
 Let us consider the steady problem (\ref{equation1}-\ref{eq1-4})
  with the local source term in the unit domain $\Om=[0,1]^2$ and set
 \ben
         \Lam= \begin{pmatrix}
        \cos\theta&\sin\theta\\
        -\sin\theta& \cos\theta
          \end{pmatrix}
          \begin{pmatrix}
        k_1&0\\
        0& k_2
          \end{pmatrix}
         \begin{pmatrix}
        \cos\theta&-\sin\theta\\
        \sin\theta& \cos\theta
         \end{pmatrix}
  \een
 and
 \ben
 f=\begin{cases}
        1000&\text{if }(x,y,z)\in[7/18,11/18]^2,\\
        0&\text{otherwise,}
    \end{cases}
 \een
 where $\theta=\pi/6$, $(k_1,k_2)=(1,1000)$.
 The boundary condition $g=0$ is imposed.

\begin{figure}[!htbp]
\centerline{\includegraphics[width=5cm,
height=5cm]{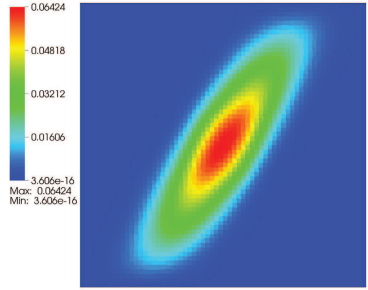}}
\caption{The numerical solution profile with monotone finite volume scheme \cite{xie}.}
\label{fig8-1}
\end{figure}
\begin{figure}[!htbp]
\begin{minipage}{0.48\linewidth}
\centerline{\includegraphics[width=6.5cm,height=5cm]{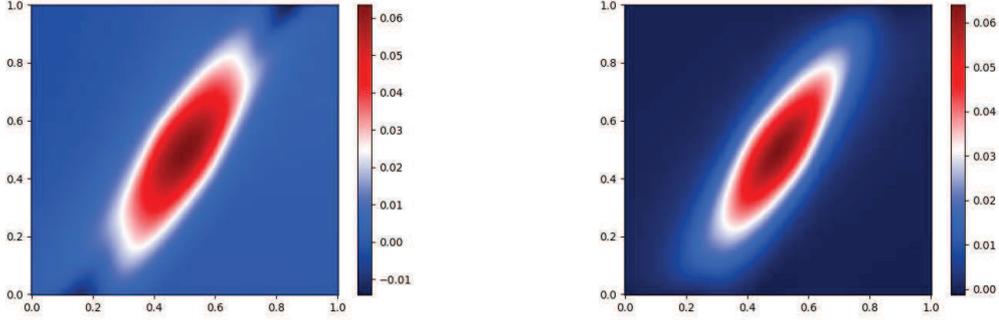}}
\end{minipage}
\hfill
\begin{minipage}{0.48\linewidth}
\centerline{\includegraphics[width=6.5cm,height=5cm]{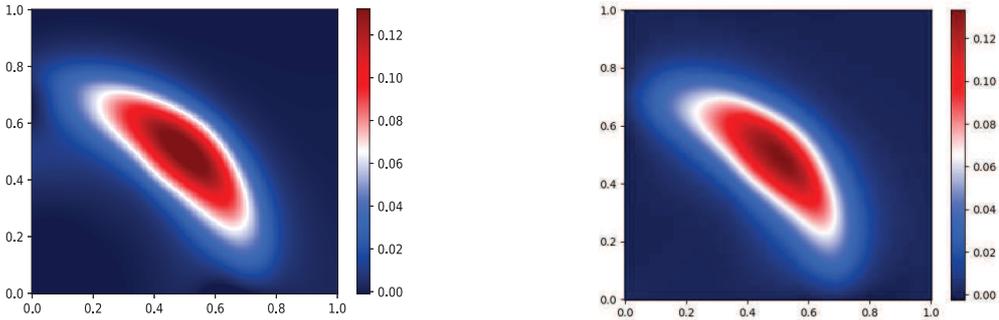}}
\end{minipage}
\caption{The solution profiles of base PINN (left) and
weighted FO-PINN (right). }
\label{fig8-2}
\end{figure}

Although there is no analytical solution, for comparison, we give the reference solution
computed with the monotone finite volume scheme \cite{xie} in
Fig. \ref{fig8-1}.
As shown in Fig. \ref{fig8-2} (left), the solution of
base PINN has more oscillations, and the heat flux is not
mainly along the rotated vector $(\sin\frac{\pi}{6},\cos\frac{\pi}{6})^T$.
while the solution of our weighted FO-PINN is much less oscillated, quite close
to the reference solution, and the heat
flux is mainly along the direction of the above rotated vector.
It is shown that our method is more accurate than the base
PINN for this rotated anisotropic diffusion problem.

Let us consider another steady problem (\ref{equation1}-\ref{eq1-4})
 with local source term in the unit domain $\Om=[0,1]^2$ and set
 \ben
         \Lam= \begin{pmatrix}
        \epsilon_1 x^2+y^2&(\epsilon_1-1)xy\\
        (\epsilon_1-1)xy & x^2+\epsilon_1 y^2
          \end{pmatrix},
          \ \epsilon_1 = 5.0e-3,
  \een
 and
 \ben
 f=\begin{cases}
        1&\text{if }(x,y)\in[3/8,5/8]^2,\\
        0&\text{otherwise.}
    \end{cases}
 \een
 The homogeneous boundary condition $g=0$ is imposed.
\begin{figure}[!htbp]
\begin{minipage}{0.48\linewidth}
\centerline{\includegraphics[width=6.5cm,height=5cm]{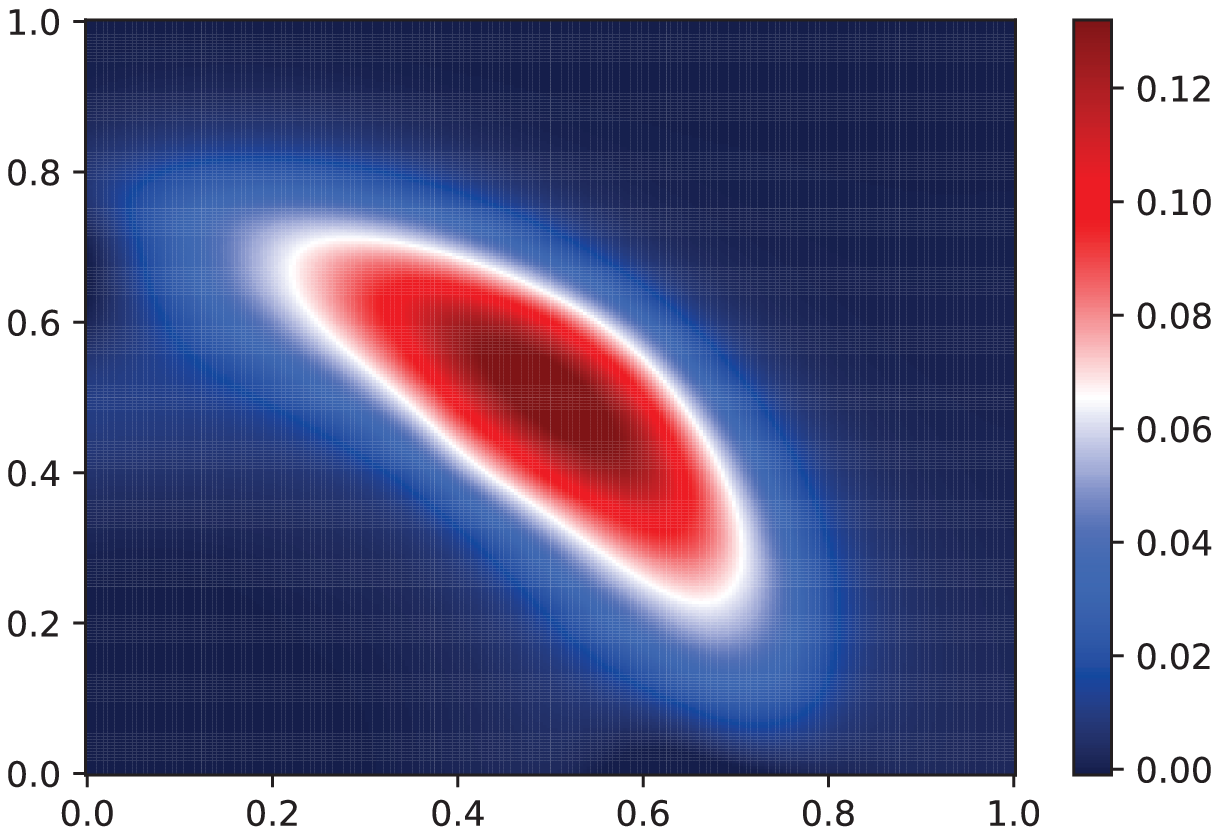}}
\end{minipage}
\hfill
\begin{minipage}{0.48\linewidth}
\centerline{\includegraphics[width=6.5cm,height=5cm]{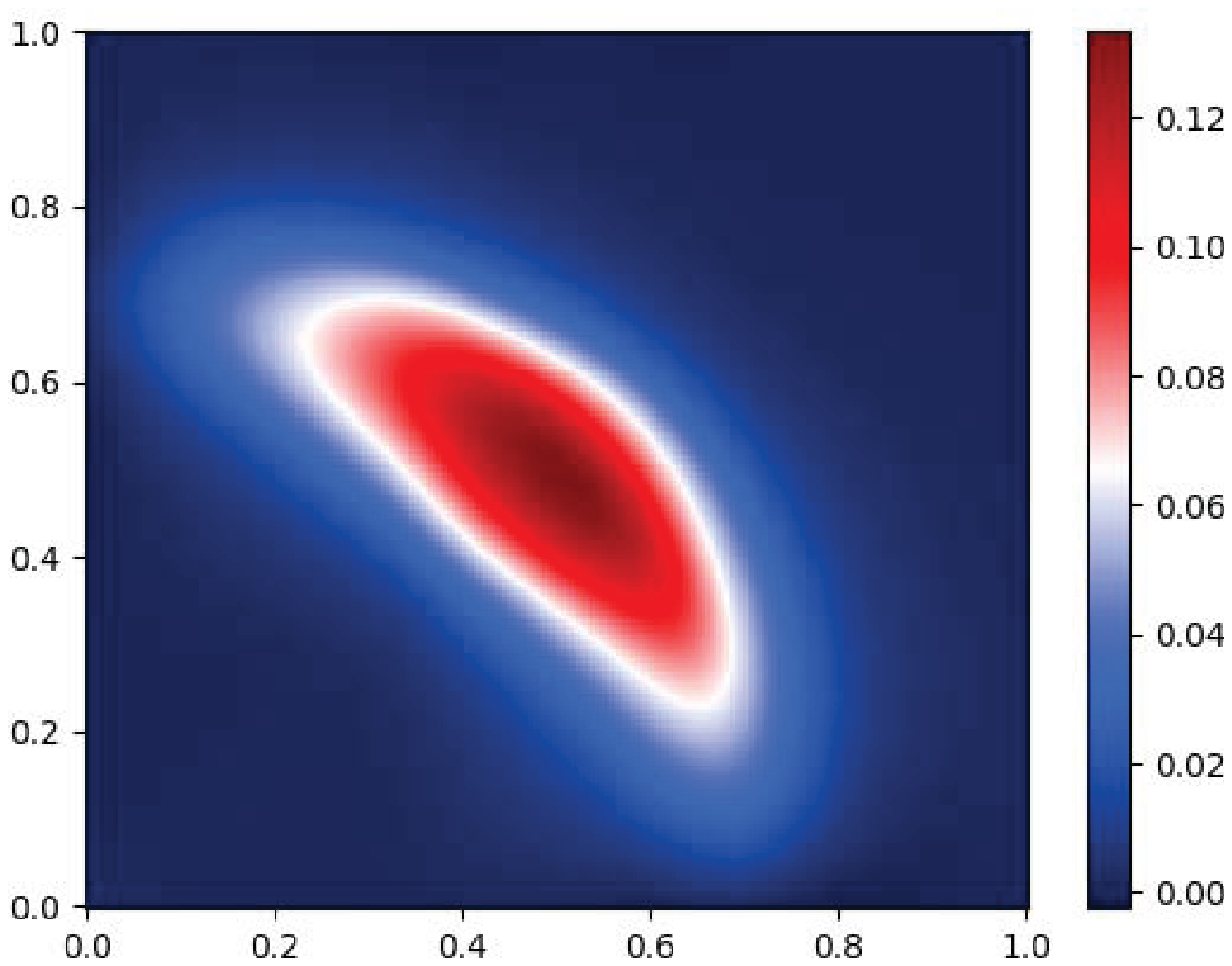}}
\end{minipage}
\caption{The solution profiles of base PINN (left) and
weighted FO-PINN (right). }
\label{fig9}
\end{figure}

With this case, we will show our method can
deal with the diffusion tensor whose
eigenvectors are
space-dependent.
As shown in Fig. \ref{fig9}, the solution
distribution of
base PINN  and our weighted FO-PINN are similar. However, compared with
the numerical results in \cite{yuan1},
the maximum solution value of our method is
more close to the one in the literature.
It is shown that our method is accurate for
this space-dependent anisotropic diffusion problem and
the accuracy is comparable to some traditional numerical schemes.

\section*{Conclusions}
A new weighted first-order formulation
is proposed for solving the anisotropic
diffusion problems with deep neural networks. The weighted
first-order formulation is first obtained by
introducing two components of the heat flux
along the eigenvectors of the
anisotropic diffusion tensor. Then the possible discontinuity of
the heat flux across the interface is correctly approximated by
multiplying the first-order formulation with the
decaying weighted function.
By various numerical examples, we can see that the smooth
diffusion problem can be solved with setting the transition width to machine
precision, and by setting a proper transition width ($\epsilon=0.01$ by
our empirical experience)
for the discontinuous interface, the numerical accuracy for
the discontinuous anisotropic diffusion problems is much improved, while the
algorithm efficiency is comparable with the base PINN due to the same
neural network structure is employed. In addition, the 3D test shows that
our method does not suffer from the curse of dimensionality as much as
other numerical schemes.
\section*{Acknowledgements}
This work was partially supported by the National Natural Science Foundation
of China (No. 11901043).


\begin{thebibliography}{17}
\bibliographystyle{plain}
\bibitem{Aavatsmark}{ I. Aavatsmark}, { An introduction to multipoint flux approximations for quadrilateral grids}, Comput. Geosci., 6 (2002), 405--432.



\bibitem{Berg}
{J. Berg, K. Nystrom}, {A unified deep artificial neural network approach to partial
differential equations in complex geometries,}
Neurocomputing, 317 (2018), 28--41.



%

\bibitem{Cai}{Z. Cai, J. Chen, M. Liu, et al.}, { Deep least-squares methods: An unsupervised learning-based numerical method for solving elliptic PDEs},
J. Comput. Phys., 420(2020), 109707.


\bibitem{chen2003}{ G. Chen, D. Li, Z. Su}, { Difference scheme by integral interpolation method
for three dimensional diffusion equations },
Chinese J. Comput. Phys., 20(2003), 205--209.

\bibitem{Cuomo}{S. Cuomo, V. S. D. Cola, F. Giampaolo, et al.}, {Scientific machine learning
through physics-informed neural networks: where we are and what's next },
(2022), arXiv: 2201.05624v3.

%
%


\bibitem{Dwivedi}
{  V. Dwivedi, B. Srinivasan}, { Physics informed extreme learning machine (PIELM)-
A rapid method for the numerical solution of partial
differential equations}, Neurocomputing, 391 (2020), 96--118.


\bibitem{Droniou1}{  J. Droniou}, {  Finite volume schemes for diffusion equations: introduction to and review of modern methods},
Math. Models Methods Appl. Sci.(M3AS), 24 (2014), 1575--1619.

\bibitem{Edwards}{ M. G. Edwards, H. Zheng}, {  Quasi M-matrix multifamily continuous Darcy-flux
approximations with full pressure support on structured and unstructured grids in three dimensions},
SIAM J. Sci. Comput., 33 (2011), 455--487.


\bibitem{Eymard}{  R. Eymard, G. Henry, R. Herbin, et al.}, {  3D bencnmark on discretization schemes for anisotropic diffusion problems on general grids},
in: Finite Volumes for Complex Applications VI Problems and Perspectives, (2011), 895--930.
%
%

\bibitem{gao2}{ Z. Gao, J. Wu}, { A linearity-preserving cell-centered scheme for the heterogeneous and
anisotropic diffusion equations on general meshes},
Int. J. Numer. Meth. Fluids, 67 (2011), 2157--2183.

%

 \bibitem{Han}{J. Han, A. Jentzen, Weinan E.},
{Solving high-dimensional partial differential equations using deep learning},
 Proc. Natl. Acad. Sci. Unit. States Am., 115(2018), 8505--8510.

 \bibitem{He}{Q. He, A. M. Tartakovsky},
{Physics-informed neural network method for forward and backward advection-dispersion
equations},
 Water Resour. Res., 57(7)(2021), e2020WR029479.

\bibitem{Jiao}{Y. Jiao, Y. Lai, D. Li, et al.}, {A rate of convergence of
physics informed neural networks for the linear
second order elliptic PDEs}, Commun. Comput. Phys., 31(4)(2022), 1272--1295.


\bibitem{Karniadakis}{G. E. Karniadakis, I. G. Kevrekidis, L. Lu, et al.}, {Physics-informed
machine learning}, Nat. Rev. Phys., 3 (2021), 422--440.

\bibitem{kershaw1981}{ D. S. Kershaw}, {  Differencing of the diffusion equation in Lagrangian hydrodynamic codes}, J. Comput. Phys., 39 (1981), 375--395.

\bibitem{Kharazmi2019}{ E. Kharazmi, Z. Zhang, G. E. Karniadakis}, { Variational physics-informed
 neural networks for solving partial differential equations},
 (2019), arXiv: 1912.00873.

\bibitem{lagaris1998}{ I. E. Lagaris, A. Likas, D. I. Fotiadis}, { Artifical neural networks for
solving ordinary and partial differential equations},
IEEE Trans. Neural Netw., 9 (1998), 987--1000.


\bibitem{Lan}{ K. Lan, J. Liu, Z. Li, et al.}, { Progress in octahedral spherical hohlraum study},
 Matter Radiat. Extrem., 1 (2016), 8--27.

  \bibitem{Lee}{H. Lee, I. S. Kang}, { Neural algorithms for solving differential equations},
J. Comput. Phys., 91 (1990), 110--131.

\bibitem{Lulu}{ L. Lu, X. Meng, Z. Mao, G. E. Karniadakis}, { DeepXDE: a deep learning
library for solving differential equations},
SIAM Rev., 63 (2021), 208--228.

%

\bibitem{Li}{ L. Li, C. Yang}, { APFOS-NET: Asymptotic preserving scheme for anisotropic
elliptic equations with deep neural network},
 J. Comput. Phys., 453 (2022), 110958.

\bibitem{Lindl}{ J. Lindl}, { Development of the indirect-drive approach to inertial confinement fusion and the target physics basis for ignition and gain},
 Phys. Plasmas, 2 (1995), 3933--4023.

\bibitem{lipnikov1}{  K. Lipnikov, M. Shashkov, D. Svyatskiy, Y. Vassilevski}, {
  Monotone finite volume schemes for diffusion equations on unstructured
 triangular and shape-regular polygonal meshes},
J. Comput. Phys., 227 (2007), 492--512.
%
%
%
%


\bibitem{Lyu}
{ L. Lyu, Z. Zhang, M. Chen, et al.}, { MIM: a deep mixed residual method for
solving high-order partial differential equations}, J. Comput. Phys., 452 (2022), 110930.

\bibitem{mao}
{Z. Mao, A. D. Jagtap, G. E. Karniadakis }, {Physics-informed neural networks for
high-speed flows}, Comput. Methods Appl. Mech. Engrg., 360 (2020), 112789.

\bibitem{meyer2012}
{C. D. Meyer, D. S. Balsara, T. D. Aslam }, {A second-order accurate Super TimeStepping
formulation for anisotropic thermal conduction}, Mon. Not. R. Astron. Soc., 422 (2012), 2102--2115.


\bibitem{Morel}{ J. Morel, J. Dendy, M. Hall, and S. White}, { A cell-centered Lagrangian-mesh diffusion differencing scheme},
J. Comput. Phys., 103 (1992), 286--299.

%
%
%


\bibitem{potier}{  C. Le Potier}, { Finite volume monotone scheme for highly anisotropic diffusion operators on unstructured triangular meshes},
C. R. Acad. Sci. Paris, Ser. I, 341 (2005), 787--792.


\bibitem{Raissi}{M. Raissi, P. Perdikaris, G. E. Karniadakis}, { Physics-informed neural networks:
a deep learning framework for solving forward and inverse problems involving
nonlinear partial differential equations},
   J. Comput. Phys., 378 (2019), 686--707.

\bibitem{schneider}{M. Schneider, B. Flemisch, R. Helmig, K. Terekhov, H. Tchelepi}, { Monotone nonlinear finite-volume method for challenging grids},
    Comput. Geosci., 22(2) (2018), 565--586.
%

\bibitem{sheng0}{ Z. Sheng, G. Yuan}, { A nine point scheme for the approximation of
    diffusion operators on distorted quadrilateral meshes},
SIAM J. Sci. Comput., 30 (2008), 1341--1361.

%

\bibitem{Shenghl}{H. Sheng, C. Yang}, {PFNN: a penalty-free neural network method for solving
a class of second-order boundary-value problems on complex geometries},
J. Comput. Phys., 428 (2021), 110085.

 \bibitem{Shin}{  Y. Shin, J. Darbon, G. E. Karniadakis}, { On the convergence of physics
 informed neural networks for linear second-order elliptic and parabolic type PDEs},
 Commun. Comput. Phys.,  28 (2020), 2042--2074.

 \bibitem{Sirignano}{J. Sirignano, K. Spiliopoulos},
{DGM: a deep learning algorithms for solving partial differential equations},
 J. Comput. Phys., 375(2018), 1339--1364.

%

\bibitem{Tartakovsky}{ A. M. Tartakovsky, C. O. Marrero, P. Perdikaris, et al.}, { Pysics-informed
deep neural networks for learning parameters and constitute relationship
in subsurface flow problems},
   Water Resour. Res., 56(5) (2020), e2019ER026731.


\bibitem{terekov}{  K. M. Terekhov, B. T. Mallison, H. A. Tchelepi}, { Cell-centered nonlinear finite volume methods for the heterogeneous anisotropic diffusion problem},
    J. Comput. Phys., 330 (2017), 245--267.

%
%

\bibitem{Wight}{ C. L. Wight, J. Zhao}, {Solving Allen-Cahn and Cahn-Hilliard equations using
the adaptive physics informed neural network}, Commun. Comput. Phys., 29 (2021), 930--954.

\bibitem{xie2}{H. Xie, C. L. Zhai, X. Xu, et al.}, { A monotone finite volume scheme with fixed stencils for 3D heat conduction equation},
Commun. Comput. Phys.,  26 (2019), 1118--1142.

\bibitem{xie}{H. Xie, X. Xu, C. L. Zhai, et al.}, { A positivity-preserving finite volume scheme
with least square interpolation for 3D anisotropic
diffusion equation},
J. Sci. Comput., 89(3),(2021), 53.



\bibitem{xiong}{F. Xiong, H. Yong, H. Chen, et al.}, { Biot's equations-based reservoir parameter
inversion using deep nerual networks}, J. Geophys. Eng., 18 (2021), 862--874.

\bibitem{yong}{H. Yong, P. Song, C. L. Zhai, et al.}, { Numerical simulation of 2-D radiation-drive
ignition implosion process}, Commun. Theor. Phys., 59 (2013), 737--744.


\bibitem{yuan1}{  G. Yuan, Z. Sheng}, {Monotone finite volume schemes for diffusion equations on polygonal meshes}, J. Comput. Phys., 227 (2008), 6288--6312.

%

\bibitem{zhao}{ F. Zhao, X. Lai, G. Yuan, Z. Sheng}, {A new interpolation for auxiliary unknowns of
    the monotone finite volume scheme for 3D diffusion equations}, Commun. Comput. Phys., 27 (2020), 1201--1233.



\end{thebibliography}
\end{document}